\documentclass[12pt]{article}
\usepackage{fullpage}
\usepackage[english]{babel}
\usepackage{amsmath,amsthm,amssymb}
\usepackage{algorithm}
\usepackage{amsfonts}
\usepackage{dsfont}
\usepackage{soul}
\usepackage{color}
\usepackage{mathrsfs}
\usepackage[flushleft]{threeparttable}
\usepackage{relsize}
\usepackage{natbib}
 \bibpunct[, ]{(}{)}{,}{a}{}{,}%
\makeatletter
\newcommand{\pushright}[1]
{\ifmeasuring@#1\else\omit\hfill$\displaystyle#1$\fi\ignorespaces}
\newcommand{\pushleft}[1]
{\ifmeasuring@#1\else\omit$\displaystyle#1$\hfill\fi\ignorespaces}
\makeatother
\allowdisplaybreaks
\newcommand{\norm}[1]{\left\lVert#1\right\rVert}
\newtheorem{thm}{Theorem}[section]

\newtheorem{prop}{Proposition}
\newtheorem{defn}{Definition}
\newtheorem{rem}{Remark}

\numberwithin{cor}{section}
\numberwithin{lem}{section}
\numberwithin{prop}{section}
\numberwithin{defn}{section}
\numberwithin{rem}{section}
\numberwithin{alg}{section}
\numberwithin{equation}{section}
\numberwithin{table}{section}

\begin{document}

\title{Solving the Dual Problems of Dynamic Programs via Regression}
\author{Helin Zhu, Fan Ye and Enlu Zhou\\
School of Industrial and Systems Engineering,\\
Georgia Institute of Technology}
\date{}
\maketitle
\begin{abstract}
In recent years, information relaxation and duality in dynamic programs have been studied extensively, and the resulted primal-dual approach has become a powerful procedure in solving dynamic programs by providing lower-upper bounds on the optimal value function. Theoretically, with the so called value-based optimal dual penalty, the optimal value function could be recovered exactly via strong duality. However, in practice, obtaining tight dual bounds usually requires good approximations of the optimal dual penalty, which could be time-consuming due to the conditional expectations that need to be estimated via nested simulation. In this paper, we will develop a framework of regression approach to approximating the optimal dual penalty in a non-nested manner, by exploring the structure of the function space consisting of all feasible dual penalties. The resulted approximations maintain to be feasible dual penalties, and thus yield valid dual bounds on the optimal value function. We show that the proposed framework is computationally efficient, and the resulted dual penalties lead to numerically tractable dual problems. Finally, we apply the framework to a high-dimensional dynamic trading problem to demonstrate its effectiveness in solving the dual problems of complex dynamic programs.\\
\\
Key words: Information relaxation, duality, dynamic program, optimal dual penalty, regression.
\end{abstract}
\section{Introduction}\label{1:intro}
Markov decision process (MDP) is a powerful model for complex dynamic decision making problems under uncertainty that arise in various fields, such as supply chain engineering, operations research, and financial engineering. The objective usually is to find an optimal policy that maximizes the expected accumulated rewards (or minimizes the expected accumulated costs) in the long run. These problems could be theoretically solved by the celebrated Bellman dynamic programming approach; however, in practice, it suffers from the so-called ``curse of dimensionality'', meaning that the size of the state space and the action space, and hence the complexity of the program increases exponentially in the dimension of the problem. Therefore, it is rarely the case that the optimal policy of a real-world dynamic programming problem could be solved exactly. Facing this issue, abundant literature has been focusing on developing good approximate dynamic programming methods that aim to construct good suboptimal policies, see \cite{bertsekas:2007}, \cite{chang:2007a} and \cite{powell:2011}, etc. In principle, given a policy (usually suboptimal), Monte Carlo simulation could be used to evaluate the policy and generate good lower bound estimators on the optimal value function by simulating a large number of state-action sample paths under the policy. It scales well with the dimension of the underlying system. However, in lack of the exact optimal value function or its upper bounds, the quality of the policy and the optimality gap of the lower bounds are difficult to measure.
\\
\\
The duality theory developed independently by \cite{rogers:2007} and \cite{brown:2010} addresses this issue by formulating and solving the dual representation of the (primal) dynamic program (DP), and providing upper bounds on the optimal value function. If the duality gap, i,e., the difference between the lower bound induced by the policy and the upper bound, is small enough, then one could claim that the policy is sufficiently good. The main idea of this duality theory is to relax the non-anticipativity constraint on all the feasible policies of the DP, i.e., allow the decision maker (DM) to choose actions based on the outcomes of future uncertainties, and penalize the DM for the access to the future information. Thus, this framework is also termed as information relaxation duality theory. In practical implementation, generating an upper bound on the optimal value function using the dual formulation only requires solving multiple pathwise deterministic inner optimization problems, and scales well with the dimension of the system via Monte Carlo simulation.
\\
\\
The information relaxation duality theory for general DPs originates from the dual theory in pricing American-style options, developed independently by \cite{rogers:2002}, \cite{haugh:2004}, and \cite{andersen:2004}, and further extended by \cite{chen2007additive}. They are able to generate upper bounds on the option price by solving the associated dual problem, which is obtained by relaxing the non-anticipativity constraint on all the feasible exercising strategies (which essentially are stopping times) and penalizing the payoff function with a martingale adapted to the nature information filtration. Furthermore, if the penalizing martingale is the Doob-Meyer martingale component of the option price process, namely the ``optimal dual martingale'', then strong duality is achieved, meaning that the upper bound is tight.
\\
\\
\cite{rogers:2007} and \cite{brown:2010} generalize the idea of duality theory from American option pricing to general discrete-time DPs and provide a broader interpretation of the dual martingale. From \cite{brown:2010}'s perspective, the dual martingale
could be regarded as the penalty for the access to the outcomes of future asset prices. Furthermore, similar to the existence of optimal dual martingale in option pricing, the authors show that in general DPs there also exist optimal dual penalties such that the dual problem with those penalties recovers the optimal value function. In particular, one specific choice of optimal dual penalties is the Doob-Meyer martingale component of the optimal value function process with respect to (w.r.t.) the natural information filtration at terminal horizon (referred to as the value-based optimal dual penalty). Following this line of research, there have been many new methodologies and applications in recent years. \cite{Brown:2011} and \cite{brown2014information} propose the gradient/subgradient penalties for the dual problems of convex DPs (concave reward function and convex action sets), in order to preserve the convexity structure in the dual problems and facilitate the optimization. They further show that these penalties could be optimal under appropriate conditions. \cite{ye2015information} generalizes the duality theory to controlled Markov diffusion (CMD) under a continuous-time setting and reveal the structure of the optimal dual penalty as a stochastic integral.
Applications include but not limited to the valuation of gas storage \cite{lai:2010}, performance measurements of trading strategies in portfolio management \cite{moallemi2015dynamic}, and robust multi-armed bandit problems \cite{kim2015robust}.
\\
\\
We also note that there are many new developments in computational methods that aim at constructing good dual penalties. In general, the (value-based) optimal dual penalty could not be computed exactly because it involves the optimal value functions that are not available and conditional expectations that need to estimated.
The naive approach is to replace the optimal value functions with  approximate ones, and use nested simulation to estimate the conditional expectations; however, this approach often requires substantial computational effort and might cause the resulted approximation to lose the dual feasibility. Various methods have been proposed to improve the accuracy and efficiency of the approximation, including the non-nested simulation approach by \cite{Belomestny:2009} and \cite{zhu2015fast} in American-style option pricing, and the pathwise optimization techniques by \cite{Desai:2011} and \cite{ye2012parameterized}. The advantage of the pathwise optimization method is that it explores a subspace of feasible dual penalties by considering the best linear combination of the existing dual penalties. However, the drawback is that it requires solving a new stochastic optimization problem, which could be computationally expensive or intractable. Moreover, the quality of the resulted dual penalty heavily relies on the quality of the existing dual penalties; hence, good performance could not be guaranteed.
\\
\\
We notice that two key things are missing in most of the existing approaches. First, the structure of the optimal dual penalty is not well-studied. Since the space consisting of all feasible dual penalties is a function space (referred to as the dual penalty space hereafter), the optimal dual penalty could be viewed as a point in that space. Therefore, we could approximate it by computing and estimating all its coordinates w.r.t. a  functional basis of the dual penalty space. If such estimation could be achieved without Monte Carlo simulation, then nested simulation is circumvented. Second, the dual problem is usually solved independently of the primal problem, meaning that some useful information (e.g., the suboptimal policy has been evaluated, and the simulations have been carried out) in the primal problem is not well-utilized. Properly utilizing that information might help facilitate estimation of the coordinates.
\\
\\
Motivated by these observations, in this paper we propose a framework of regression approach that explores the structure of the dual penalty space as well as the optimal dual penalty while utilizing the information from the primal problem. In particular, we find two function bases of the dual penalty space that possess desired properties such that the corresponding coordinates of the optimal dual penalty are easy to compute and the resulted approximations of the optimal dual penalty remain dual feasible. Furthermore, instead of using Monte Carlo simulation to estimate those coordinates, we propose an efficient regression method that reuses the suboptimal policy and the samples generated in the primal problem.
\\
\\
The proposed framework of regression approach has the following advantages. First, it is a scheme that requires minimal extra simulation or computational costs since it reuses the suboptimal policy and samples generated in the primal problem. Second, it generates approximations of the optimal dual penalty without nested simulation suffered by some existing approaches. Therefore, computational efficiency is significantly improved. Third, in contrast to some of the existing approaches, the resulted dual penalties are feasible, and thus induce valid upper bounds on the optimal value function. Finally, as we will see later, the framework is robust. Dual penalties with desired structural properties could be generated by constructing proper regressors in the regression. For instance, for a convex DP, the framework is able to generate dual penalties that preserve convexity in the dual problem.
\\
\\
It turns out, in a broader sense, several existing approaches to approximating the optimal dual penalty could be regarded as special cases of this framework under specific settings. For instance,  \cite{Belomestny:2009} proposes a non-nested simulation approach to approximating the optimal dual martingale in American option pricing. In particular, the authors apply martingale representation theorem on the optimal dual martingale, then approximate the resulted stochastic integral with an Ito sum, and finally estimate the integrands in the Ito sum via regression. We will show that this approach could be viewed as a special case of the proposed framework with a specific functional basis of the dual penalty space. Our framework is more universal and powerful, because it reveals the structure of the optimal dual penalty regardless of the underlying probability measure (i.e., not restricted to the Brownian measure in \cite{Belomestny:2009} or the Poisson random measure in \cite{zhu2015fast}). We will also show that the approximation scheme in \cite{ye2015information} could be viewed as a special case of the proposed framework as well. To summarize, the contributions are as follows:
\begin{itemize}
\item We develop a framework of regression approach to approximating the optimal dual penalty in general DPs that (1) circumvents nested simulation that is suffered by some existing approaches, (2) requires minimal extra simulation or computational costs by reusing the suboptimal policy and samples generated in the primal problem, and (3) yields feasible dual penalties and valid upper bounds. Furthermore, the framework is quite robust in the sense that it is capable of generating dual penalties with desired structural properties by properly constructing the regressors in the regression.

\item We explore the dual penalty space and reveal the structure of the optimal dual penalty for general DPs, which enable us to generate good feasible dual penalties in a systematic way. Several existing approaches could be regarded as special cases of the framework under specific settings.

\item The application of the proposed framework to a high-dimensional dynamic trading problem shows that it is effective and efficient in solving the dual problems of complex DPs. In particular, it generates accurate approximations of the optimal dual penalty and tight upper bounds on the optimal value function when good suboptimal policies and appropriate regressors are used in the regression approach.
\end{itemize}
The rest of the paper is organized as follows. In Section \ref{sec2:Formulation}, we review the basics of dynamic programming and information relaxation duality theory. We present the framework of regression approach in Section \ref{sec3:regression}. In Section \ref{sec4:special}, we show several existing approaches are special cases of the proposed framework under specific settings. In Section \ref{sec5:numerical}, we apply the proposed framework to a high-dimensional dynamic trading problem, and demonstrate its effectiveness and efficiency in solving the dual problems of complex DPs. Conclusions are provided in Section \ref{sec6:conclusion}.
\section{Dynamic Programs and Dual Formulations}\label{sec2:Formulation}
On a general probability space $(\Omega,\mathbb{P},\mathscr{F})$, in which $\Omega$ is the set of all possible outcomes (scenarios) of uncertainties, $\mathbb{P}$ is the underlying probability measure and $\mathscr{F}$ is the $\sigma$-algebra consisting of all the events (measurable subsets of $\Omega$). Consider a finite-horizon MDP as follows. Time is indexed by $\mathcal{T}=\{0,1,...,N\}$. The state $x$ follows the dynamics
\begin{equation}\label{eq.2.1}
x_{n+1}=f(x_n, a_n, z_{n+1}), \quad n=0,1,\cdots,N-1,
\end{equation}
where $f$ is the deterministic transition function, $x_{n}\in \mathscr{X}_{n}$ denotes the state at period $n$ that lives in the state space $\mathscr{X}_{n}$, $a_{n}\in \mathscr{A}_{n}$ denotes the action/control at period $n$ that is chosen from the action space $\mathscr{A}_{n}$, $z_{n+1}$ is the random noise at period $n$ and $\{z_{n}:n=1,...,N\}$ are assumed to be independently and identically distributed (i.i.d.) random variables with probability measure $\rho$ on support space $\Xi\subseteq\mathbb{R}^d$. The evolution of information is described by the natural filtration $\mathcal{F}=\{\mathscr{F}_{n}: n=0,...,N\}$. Loosely speaking, $\mathscr{F}_{n}$ describes the information available to the DM at period $n$. In particular, each $z_{n}$ is $\mathscr{F}_n$-measurable. Without loss of generality, we further assume $\mathscr{F}_{0}=\{\emptyset,\Omega\}$, meaning that the DM initially does not have any knowledge about the outcomes of uncertainties, and $\mathscr{F}_{N}=\mathscr{F}$, meaning that the DM knows all the possible outcomes of uncertainties in the end.
\\
\\
We use a mapping $\alpha_{n}(\cdot)$ from the state space $\mathscr{X}_{n}$ to the action space $\mathscr{A}_{n}$, i.e., $\alpha_{n}:\mathscr{X}_{n}\rightarrow \mathscr{A}_{n}$, to denote the decision rule at period $n$. A policy/strategy $\boldsymbol{\alpha}:=(\alpha_{0},\alpha_1,...,\alpha_{N-1})$ consisting of a sequence of decision rules is called non-anticipative/$\mathcal{F}$-adapted, if each decision rule $\alpha_{n}(\cdot)$ is $\mathscr{F}_{n}$-measurable. Intuitively, it means the DM chooses the action $a_{n}$ only based on the information accumulated up to period $n$, he/she shall not choose the action $a_{n}$ based on the future information. We use $\mathbb{A}_{\mathcal{F}}$ to denote the set of all non-anticipative policies, and $\mathbb{A}$ to denote the set of all policies (including the anticipative ones); clearly $\mathbb{A}_{\mathcal{F}}\subseteq \mathbb{A}$. Furthermore, we associate an $\mathscr{F}_{n}$-measurable reward function $r_{n}(x_{n},a_{n})$ with the state dynamics to represent the immediate reward after the DM chooses action $a_{n}$ at period $n=0,...N-1$, and an $\mathscr{F}_{N}$-measurable function $r_{N}(x_{N})$ as the terminal reward at period $N$.
\\
\\
Given $x_0\in \mathscr{X}_0$, the objective of the DM is to select a non-anticipative policy $\boldsymbol{\alpha}\in \mathbb{A}_{\mathcal{F}}$ that maximizes the accumulated rewards over all the periods, i.e.,
\begin{equation}\label{eq.2.2}
(P): V_{0}(x_{0})\overset{\triangle}=
\sup\limits_{\mathbf{\boldsymbol\alpha}\in \mathbb{A}_{\mathcal{F}}}\mathbb{E}_{0}
\left[\sum\limits_{n=0}^{N-1}r_n(x_{n},\alpha_{n}(x_{n}))
+r_{N}(x_{N})\right],
\end{equation}
where $\mathbb{E}_{0}\left[\cdot\right]$ means that the expectation is taken w.r.t. $\mathscr{F}_{0}$, and we will use $\mathbb{E}_{n}\left[\cdot\right]$ to denote the expectation taken w.r.t. $\mathscr{F}_{n}$ thereafter.
\\
\\
It is well known that problem (\ref{eq.2.2}) could be recursively solved theoretically via Bellman backward dynamic programming
\begin{equation}\label{eq.2.3}
\left\{\begin{array}{l}
V_{N}\left(x_{N}\right)\overset{\triangle}=r_{N}(x_{N}),\\
V_{n}\left(x_{n}\right)\overset{\triangle}=\sup\limits_{a_{n}\in \mathscr{A}_{n}}\left\{r_{n}\left(x_{n},a_{n}\right)+
\mathbb{E}_{n}\left[V_{n+1}\left(x_{n+1}\right)\right]\right\}, \;n=N-1,...,0,
\end{array}\right.
\end{equation}
where $V_{n}\left(x_{n}\right)$ represents the optimal value function of the DP with initial period $n$ and initial state $x_{n}$. However, in practice, the Bellman recursion (\ref{eq.2.3}) could hardly be solved exactly for most cases, due to the curse of dimensionality. Therefore, one often needs to settle with suboptimal policies or approximate optimal value functions. A suboptimal policy could be evaluated conveniently via Monte Carlo simulation to generate a lower bound on the optimal value function. However, in the absence of the exact optimal value function or its upper bounds, the quality of the suboptimal policy could hardly be measured.
\subsection{A Dual Formulation}\label{sec2.a: dual formulation}
The duality theory developed by \cite{rogers:2007} and \cite{brown:2010} addresses the aforementioned issue by formulating a dual representation of the (primal) DP (\ref{eq.2.2}), and solving the dual problem provides a valid upper bound on the optimal value function. Therefore, the quality of a suboptimal policy could be empirically measured by examining the duality gap. If it is sufficiently small, then the suboptimal policy could be claimed to be near optimal.
\\
\\
To be more specific, let us rigorously define a \emph{feasible dual penalty} as follows.
\begin{defn}\label{def.2.1}
We say $M(\boldsymbol{\alpha},\mathbf{z})$, a functional of policy $\boldsymbol\alpha \in \mathbb{A}$ and noise sequence $\mathbf{z}:= (z_{1},...,z_{N})$, is a feasible dual penalty if
\begin{equation}\label{eq.2.4}
\mathbb{E}_{0}\left[M(\boldsymbol{\alpha},\mathbf{z})\right] = 0, \quad \forall \boldsymbol{\alpha}\in \mathbb{A}_{\mathcal{F}}.
\end{equation}
\end{defn}
\noindent Put in another way, a penalty function $M(\boldsymbol{\alpha},\mathbf{z})$ is dual feasible if it does not penalize any non-anticipative policy in expectation. We further use $\mathbb{M}_{\mathcal{F}}$ to denote the set of all feasible dual penalties.
\begin{rem}\label{rem.2.1}
Definition \ref{def.2.1} is slightly different from the one in \cite{brown:2010}, in which a dual penalty is called feasible if
$\mathbb{E}_{0}\left[M(\boldsymbol{\alpha},\mathbf{z})\right]\le 0, \; \forall \boldsymbol{\alpha}\in \mathbb{A}_{\mathcal{F}}$. The reason for using Definition \ref{def.2.1} is to ensure that the set of all feasible dual penalties $\mathbb{M}_{\mathcal{F}}$ is a function space (vector space). Note that Definition \ref{def.2.1} does not exclude any ``good'' feasible dual penalties in \cite{brown:2010}. If $M(\boldsymbol{\alpha},\mathbf{z})$ is a feasible dual penalty in \cite{brown:2010}, i.e., $\mathbb{E}_{0}\left[M(\boldsymbol{\alpha},\mathbf{z})\right]\le 0, \; \forall \boldsymbol{\alpha}\in \mathbb{A}_{\mathcal{F}}$. Then $\mathbb{E}_{0}\left[(M(\boldsymbol{\alpha},\mathbf{z})-
\mathbb{E}_{0}\left[M(\boldsymbol{\alpha},\mathbf{z})\right])\right]=0, \; \forall \boldsymbol{\alpha}\in \mathbb{A}_{\mathcal{F}}$. That is, $(M(\boldsymbol{\alpha},\mathbf{z})-
\mathbb{E}_{0}\left[M(\boldsymbol{\alpha},\mathbf{z})\right])$ is a feasible dual penalty by Definition \ref{def.2.1}, and it always induces an upper bound as tight as the one induced by $M(\boldsymbol{\alpha},\mathbf{z})$. In particular, Definition \ref{def.2.1} does not exclude the optimal dual penalty (which will be defined in the following).
\end{rem}

\noindent \cite{brown:2010} shows that the following dual representation of the primal DP (\ref{eq.2.2}), derived by subtracting a feasible dual penalty $M(\boldsymbol\alpha,\mathbf{z})\in\mathbb{M}_{\mathcal{F}}$ from the objective function, relax the non-anticipativity constraint on all the feasible policies, and interchange the maximization and expectation, i.e.,
\begin{equation}\label{eq.2.5}
(D): V_{0}^{M}(x_{0})
\overset{\triangle}=
\mathbb{E}_{0}\left[\sup_{\mathbf{a}\in \mathscr{A}} \left\{\sum\limits_{n=0}^{N-1}
r_{n}\left(x_{n},a_{n}\right)+r_{N}(x_{N})
-M(\boldsymbol{a},\mathbf{z})\right\}\right],
\end{equation}
yields an upper bound $V_{0}^{M}(x_{0})$ on the optimal value function $V_{0}(x_{0})$, where $\mathbf{a}:=(a_0,...,a_{N-1})$ represents the action sequence and $\mathscr{A}:=(\mathscr{A}_0,...,\mathscr{A}_{N-1})$. Conceptually, the dual problem (\ref{eq.2.5}) consists of a series of scenario-based deterministic pathwise optimization problems (referred to as inner optimization problems), in which the actions could be chosen with full knowledge of all outcomes of uncertainties.
In practice, solving dual problem (\ref{eq.2.5}) is convenient via Monte Carlo simulation given that the inner optimization problems are tractable: simulate multiple i.i.d. noise sequences $\mathbf{z}$, then solve the deterministic inner optimization problem corresponding to each $\mathbf{z}$, and finally take the average of the optimal values as an estimator of $V_{0}^{M}(x_{0})$.
\\
\\
Note that a dual penalty is feasible if and only if it does not penalize any non-anticipative policy, we immediately have the weak duality. Furthermore, if minimizing $V_{0}^{M}(x_{0})$ over all feasible dual penalties $M\in \mathbb{M}_{\mathcal{F}}$, then we have $\forall x_0\in \mathscr{X}_0$,
$$
\sup\limits_{\mathbf{\boldsymbol\alpha}\in \mathbb{A}_{\mathcal{F}}}V_{0}^{\boldsymbol{\alpha}}(x_{0})=
V_{0}(x_{0})=\inf\limits_{M \in \mathbb{M}_{\mathcal{F}}}
V_{0}^{M}(x_{0}),
$$
i.e., strong duality holds. In particular, the following proposition shows that the value-based optimal dual penalty in the form of (\ref{eq.2.6}) suffices for the strong duality to hold.
\begin{prop}\label{prop.2.1}
(i) (\textbf{Weak Duality})[Lemma 2.1 in \cite{brown:2010}]
If $M(\boldsymbol{\alpha},\mathbf{z}) \in \mathbb{M}_{\mathcal{F}}$, then we have
$$
V_{0}(x_{0})\le V_{0}^{M}(x_{0}).
$$
(ii) (\textbf{Strong Duality}) [Theorem 2.3 in \cite{brown:2010}]
Let $\{V_{n}(x_{n})\}$ be the optimal value functions of the primal DP (\ref{eq.2.2}). Further let $M^{\ast}(\boldsymbol{\alpha},\mathbf{z})$ be the martingale difference sum of $\{V_{n}(x_{n})\}$, i.e.,
\begin{equation}\label{eq.2.6}
M^{\ast}(\boldsymbol{\alpha},\mathbf{z})\overset{\triangle}=
\sum\limits_{n=0}^{N-1}\left(V_{n+1}(x_{n+1})-
\mathbb{E}\left[V_{n+1}(x_{n+1})|x_n, a_n\right]\right),\quad
\text{where}\; a_n=\alpha_n(x_n),
\end{equation}
then $M^{\ast}(\boldsymbol{\alpha},\mathbf{z})$ is a feasible dual penalty. Moreover, it achieves strong duality, i.e.,
\begin{equation}\label{eq.2.7}
V_{0}^{M^{\ast}}(x_{0})
\overset{\triangle}=\mathbb{E}_{0}\left[ \sup\limits_{\boldsymbol\alpha\in \mathbb{A}}\left\{\sum\limits_{n=0}^{N-1}r_{n}
\left(x_{n},\alpha_{n}(x_{n})\right)+
r_{N}(x_{N})-M^{\ast}(\boldsymbol{\alpha},\mathbf{z})\right\} \right]=V_{0}(x_{0}).
\end{equation}
\end{prop}
\noindent Note that the optimal dual penalty $M^{\ast}(\boldsymbol{\alpha},\mathbf{z})$ defined in (\ref{eq.2.6}) is a functional of
$(\boldsymbol{\alpha},\mathbf{z})$ because the state sequence $(x_0,...,x_N)$ and action sequence $\mathbf{a}$ depend on $(\boldsymbol{\alpha},\mathbf{z})$ through the state dynamics \eqref{eq.2.1} and the policy. Intuitively, weak duality suggests that the dual problem \eqref{eq.2.5} with any feasible dual penalty could be used to generate an upper bound on the optimal value function, and strong duality implies that the advantage gained by the access to the future information is perfectly cancelled out in expectation by the optimal dual penalty. What is more striking about the optimal dual penalty is that the second equality in (\ref{eq.2.7}) is achieved almost surely for every inner optimization problem scenario. Therefore, we can drop the expectation sign in (\ref{eq.2.7}) and strengthen the result as
\begin{equation*}
V_{0}(x_{0})=\sup\limits_{\boldsymbol\alpha\in \mathbb{A}}\;
\biggl\{\sum\limits_{n=0}^{N-1}r_{n}\left(x_{n},\alpha_{n}(x_{n})\right)
+r_{N}(x_{N})
-M^{\ast}(\boldsymbol{\alpha},\mathbf{z})\biggr\},\quad
a.s..
\end{equation*}
Therefore, in practical implementation, if the approximation of the optimal dual penalty is sufficiently accurate, then the variance of a one-sample estimator of the upper bound is small. Therefore, one only needs to solve a small number of deterministic inner optimization problems to generate a good upper bound estimator.
\\
\\
In theory, strong duality implies that the upper bound induced by the optimal dual penalty is tight. However, one could hardly compute the optimal dual penalty exactly and close the duality gap since in general the optimal value functions $\{V_{n}(x_{n})\}$ are not available. A naive alternative is to replace the optimal value functions with approximate ones $\{\widetilde{V}_{n}(x_{n})\}$ that might be induced by a suboptimal policy or certain approximate dynamic programming technique. The computational effort required in approximating the optimal dual penalty should be taken into consideration as well. A common method to estimate the conditional expectations in the optimal dual penalty \eqref{eq.2.6} is nested simulation, which generates scenarios in the outer-layer and uses sample averaging in the inner-layer associated with each scenario. Hence, the total simulation effort is proportional to the number of outer-layer scenarios multiplied by the number of inner-layer samples. Albeit stable with a large-scale simulation, it might not be a desirable choice facing a limited computational budget. Several approaches have been developed in recent years that aim to generate good dual penalties without nested simulation. However, we note that two key things are missing in most of the existing approaches: the structure of the dual penalty space as well as the optimal dual penalty is not well-studied; the information such as the suboptimal policy and the sample paths generated in solving the primal problem is not well-utilized.
\\
\\
In the next section, we will present a framework of regression approach to approximating the optimal dual penalty that explores the structure of dual penalty space and efficiently reutilizes the information in the primal problem.
\section{A Framework of Regression Approach}\label{sec3:regression}
\subsection{General Framework}\label{sec3.a:framework}
The main idea of the regression approach is to view the optimal dual penalty as a point in the dual penalty space, then compute its coordinates w.r.t. a properly chosen functional basis, and finally estimate those coordinates in a non-nested manner. The key question lies in how to choose the functional basis and compute the corresponding coordinates of the optimal dual penalty. Let us first derive the general framework. Recall the optimal dual penalty
\begin{equation*}
M^{\ast}(\boldsymbol{\alpha},\mathbf{z})=
\sum\limits_{n=0}^{N-1}\left(V_{n+1}(x_{n+1})-
\mathbb{E}\left[V_{n+1}(x_{n+1})|x_n, a_n\right]\right),\quad \text{where}\; a_n=\alpha_n(x_n).
\end{equation*}
Denote the $n$-th single martingale difference term in $M^{\ast}(\boldsymbol{\alpha},\mathbf{z})$ by $M^{\ast}_n(\boldsymbol{\alpha},\mathbf{z})$, i.e.,
\begin{equation}\label{eq.3.1.2}
M^{\ast}_n(\boldsymbol{\alpha},\mathbf{z})\overset{\triangle}=
\left(V_{n+1}(x_{n+1})-
\mathbb{E}\left[V_{n+1}(x_{n+1})|x_n, a_n\right]\right).
\end{equation}
Then clearly $M^{\ast}=\sum_{n=0}^{N-1}M^{\ast}_n$. Further let $\mathcal{B}=\{b_i: i\in I\}$ denote a functional basis of $\mathbb{M}_{\mathcal{F}}$, where $\{b_i(\cdot)\}$ are functions on the support space $\Xi$ of the noises $\{z_n\}$ such that $\int b_i(z)\rho(dz)=0, \forall i\in I$, where recall that $\rho$ is the probability measure induced by the noises $\{z_n\}$. Therefore, we have
\begin{equation}\label{eq.3.3}
\mathbb{E}\left[b_i(z_{n+1})\right]=0, \; \forall i\in I, n=0,...,N-1.
\end{equation}
Note that $\mathcal{B}$'s cardinality, $|I|$, might be countable or uncountable. We will discuss about the construction of $\mathcal{B}$ or the properties it should possess later. Suppose we could express $M^{\ast}_n(\boldsymbol{\alpha},\mathbf{z})$ w.r.t. $\mathcal{B}$ in the form of
\begin{equation}\label{eq.3.1.3}
\left(V_{n+1}(x_{n+1})-
\mathbb{E}\left[V_{n+1}(x_{n+1})|x_n, a_n\right]\right)
=\sum_{i \in I}\beta_{n,i}(x_n, a_n)\cdot b_i(z_{n+1}),
\end{equation}
where $\{\beta_{n,i}(x_n, a_n)\}$ are the coordinates. To compute $\{\beta_{n,i}(x_n, a_n)\}$, multiplying $b_i(z_{n+1})$ $\forall i\in I$ on both sides of (\ref{eq.3.1.3}) and taking conditional expectations w.r.t. $\mathscr{F}_n$, we obtain
\begin{equation}\label{eq.3.4}
\begin{split}
&\mathbb{E}_n\left[V_{n+1}(x_{n+1})\cdot b_i(z_{n+1})\right]
-\mathbb{E}_n\left[\mathbb{E}\left[V_{n+1}(x_{n+1})|x_n, a_n\right]\cdot b_i(z_{n+1})\right]\\
={}&\sum_{j \in I}\mathbb{E}_n\left[\beta_{n,j}(x_n, a_n)\cdot
b_i(z_{n+1})\cdot
b_j(z_{n+1})\right].
\end{split}
\end{equation}
Due to the Markov property, we have
\begin{equation*}
\mathbb{E}_n\left[V_{n+1}(x_{n+1})\cdot b_i(z_{n+1})\right]
=\mathbb{E}\left[V_{n+1}(x_{n+1})\cdot b_i(z_{n+1})|x_n, a_n\right].
\end{equation*}
Furthermore, by (\ref{eq.3.3}) and the facts that $\mathbb{E}\left[V_{n+1}(x_{n+1})|x_n, a_n\right]$ is $\mathscr{F}_n$-measurable and $z_{n+1}$ is independent of $\mathscr{F}_n$, we have
\begin{equation*}
\begin{split}
&\mathbb{E}_n\left[\mathbb{E}\left[V_{n+1}(x_{n+1})|x_n, a_n\right]\cdot b_i(z_{n+1})\right]
=\mathbb{E}\left[V_{n+1}(x_{n+1})|x_n, a_n\right]
\cdot \mathbb{E}_n\left[b_i(z_{n+1})\right]\\
={}& \mathbb{E}\left[V_{n+1}(x_{n+1})|x_n, a_n\right]
\cdot \mathbb{E}\left[b_i(z_{n+1})\right]
=0.
\end{split}
\end{equation*}
Hence, equation (\ref{eq.3.4}) could be rewritten as
\begin{equation}\label{eq.3.1.4}
\mathbb{E}\left[V_{n+1}(x_{n+1})\cdot b_i(z_{n+1})|x_n, a_n\right]
=\sum_{j \in I}\beta_{n,j}(x_n, a_n)\mathbb{E}\left[b_i(z_{n+1})\cdot
b_j(z_{n+1})\right],
\end{equation}
where we use the facts that $\beta_{n,j}(x_n, a_n)$ is $\mathscr{F}_n$-measurable and $z_{n+1}$ is independent of $\mathscr{F}_n$. Note that $\{\mathbb{E}\left[b_i(z_{n+1})\cdot
b_j(z_{n+1})\right]\}$ in (\ref{eq.3.1.4}) are constants that only depend on $\mathcal{B}$. Therefore, in principle we could view (\ref{eq.3.1.4}) as a system of linear equations with variables $\{\beta_{n,i}(x_n, a_n)\}$ (although the number of equations in the system, $|I|$, could be countably infinite or even uncountable). Assuming it is solvable (again this depends on the choice of $\mathcal{B}$), then we have
\begin{equation}\label{eq.3.1.5}
\beta_{n,i}(x_n, a_n)=\mathbb{E}
\left[V_{n+1}(x_{n+1})\cdot h_i(z_{n+1})|x_n, a_n\right],
\quad \forall i\in I,\; n=0,...,N-1,
\end{equation}
where $\{h_i(\cdot): i\in I\}$ are deterministic functions that only depend on $\mathcal{B}$ because all the parameters of the linear system are uniquely determined by $\mathcal{B}$. It follows that the optimal dual penalty could be rewritten as
\begin{equation}\label{eq.3.1.6}
M^{\ast}(\boldsymbol{\alpha},\mathbf{z})=\sum_{n=0}^{N-1}
\sum_{i\in I}\beta_{n,i}(x_n, a_n)\cdot b_i(z_{n+1}),
\quad \text{where}\; a_n=\alpha_n(x_n),
\end{equation}
where the coordinates $\{\beta_{n,i}(x_n, a_n)\}$ are given by (\ref{eq.3.1.5}).
\\
\\
Instead of approximating the optimal dual penalty $M^{\ast}(\boldsymbol{\alpha},\mathbf{z})$ directly via nested simulation, now we are able to approximate it by estimating all its coordinates $\{\beta_{n,i}(x_n, a_n)\}$ w.r.t. $\mathcal{B}$. The difficulty lies in how to efficiently estimate the conditional expectations on the right hand side of (\ref{eq.3.1.5}). To avoid nested simulation, we propose the following regression scheme. The main idea is to treat $\beta_{n,i}(x_n, a_n)$ as the expected response, and thus in the regression an observation of the expected response is a sample outcome of $V_{n+1}(x_{n+1})\cdot h_i(z_{n+1})|x_n, a_n$, which is obtained by exercising a given suboptimal policy $\widetilde{\boldsymbol{\alpha}}:=
\left(\widetilde{\alpha}_0,\widetilde{\alpha}_1,...,
\widetilde{\alpha}_{N-1}\right)$ and computing the value function along one sample path.
Further let
$\boldsymbol{\phi}_{n,i}(x_n, a_n)\overset{\triangle}=(\phi_{n,i}^1(x_n, a_n),..., \phi_{n,i}^K(x_n, a_n))^T$ denote the vector of regressors (dependent variables) in the regression, and later we will illustrate how to properly construct them based on the choice of $\mathcal{B}$. The linear regression model could be formulated as
\begin{equation*}
V_{n+1}(x_{n+1})\cdot h_i(z_{n+1})|x_n, a_n=
\left(\boldsymbol{\phi}_{n,i}(x_n, a_n)\right)^T \boldsymbol{\theta}_{n,i}+\epsilon_{n,i},
\end{equation*}
where $\boldsymbol{\theta}_{n,i}:=
(\theta_{n,i}^1,...,\theta_{n,i}^K)^T$ is the vector of regression coefficients and $\epsilon_{n,i}$ is the noise.
The complete algorithm is summarized as the following Algorithm \ref{alg.3.1}.

\begin{algorithm}
\caption{A Regression-based Algorithm for Approximating the Optimal Dual Penalty}\label{alg.3.1}
\textbf{Input}: Functional basis $\mathcal{B}$ and suboptimal policy $\widetilde{\boldsymbol{\alpha}}$.\\
\textbf{Output}: $\widetilde{M}(\boldsymbol{\alpha},\mathbf{z})$---Approximation of the optimal dual penalty.\\
\textbf{1. Initialization}: Simulate $M$ independent (state-action) sample paths under the given policy $\widetilde{\boldsymbol{\alpha}}$ subject to state dynamics (\ref{eq.2.1}). Denote the sample paths by
$$\{(x_0^j,a_0^j;x_1^j,a_1^j;...,x_{N-1}^j,a_{N-1}^j; x_{N}^j): j=1,...,M\},$$
where the noise sequences are $\{(z_1^j,z_2^j,...,z_N^j): j=1,...,M\}$ and $a_n^j=\widetilde{\alpha}_n(x_n^j)$. Calculate the corresponding value functions along each sample path, denoted by
$$\{(V_0^j(x_0^j),V_1^j(x_1^j),...,V_N^j(x_N^j)): j=1,...,M\},$$ where $
V_n^j(x_n^j)\overset{\triangle}
=\sum_{k=n}^{N-1}r_n(x_n^j,a_n^j)+r_N(x_N^j).
$
\\
\textbf{2. Iteration}: For $n=0,...,N-1$, $i\in I$, estimate the coordinate $\beta_{n,i}(x_n, a_n)$ by regression. The responses are $\{V_{n+1}^j(x_{n+1}^j)\cdot h_i(z_{n+1}^j):j=1,...,M\}$ and the design matrix is
$(\boldsymbol{\phi}^T_{n,i}(x^1_n, a^1_n),...,
\boldsymbol{\phi}^T_{n,i}(x^M_n, a^M_n))^T$,
where $\boldsymbol{\phi}_{n,i}(x^j_n, a^j_n)=
(\phi_{n,i}^1(x^j_n, a^j_n),..., \phi_{n,i}^K(x^j_n, a^j_n))^T$. Let $\widetilde{\boldsymbol{\theta}}_{n,i}:=
(\widetilde{\theta}_{n,i}^1,...,\widetilde{\theta}_{n,i}^K)^T$ denote the vector of regression coefficients,
i.e.,
\begin{equation*}
\widetilde{\boldsymbol{\theta}}_{n,i}=
    \arg\min_{\boldsymbol{\theta}_{n,i}}\frac{1}{M}\sum_{j=1}^M
    \biggl(V_{n+1}^j(x_{n+1}^j)\cdot h_i(z_{n+1}^j)-\left(\boldsymbol{\phi}_{n,i}(x^j_n, a^j_n)\right)^T\boldsymbol{\theta}_{n,i} \biggr)^2.
\end{equation*}
\textbf{3. Termination}: Let
\begin{equation*}
\widetilde{\beta}_{n,i}(x_n,a_n)\overset{\triangle}
=\left(\boldsymbol{\phi}_{n,i}(x_n, a_n)\right)^T\widetilde{\boldsymbol{\theta}}_{n,i} .
\end{equation*}
be the estimation of $\beta_{n,i}(x_n,a_n)$ and
\begin{equation}\label{eq.3.1.9}
\widetilde{M}(\boldsymbol{\alpha},\mathbf{z})\overset{\triangle}=
\sum_{n=0}^{N-1}\sum_{i\in I} \widetilde{\beta}_{n,i}(x_n,
a_n)\cdot b_i(z_{n+1}), \quad \mbox{where}\; a_n=\alpha_n(x_n)
\end{equation}
be the approximation of $M^{\ast}(\boldsymbol{\alpha},\mathbf{z})$.
\end{algorithm}

\noindent
\\
For convenience, we refer to $\widetilde{M}(\boldsymbol{\alpha},\mathbf{z})$ as the \emph{regression-based (dual) penalty}. Note that one could implement Algorithm \ref{alg.3.1} efficiently by reusing the sample paths that were generated to evaluate the suboptimal policy $\widetilde{\boldsymbol{\alpha}}$ in the primal problem. Hence,  minimal extra simulation or computational costs are required. Moreover, generating a regression-based dual penalty $\widetilde{M}(\boldsymbol{\alpha},\mathbf{z})$ as in (\ref{eq.3.1.9}) does not incur nested simulation, since the conditional expectations are estimated via regression instead. Hence, computational efficiency is significantly improved. Finally, one of the biggest advantages of the proposed regression approach over some of the existing approaches is that it preserves the dual feasibility of the optimal dual penalty as long as the function basis $\mathcal{B}$ satisfies (\ref{eq.3.3}). That is, the regression-based penalty $\widetilde{M}(\boldsymbol{\alpha},\mathbf{z})$ maintains to be a feasible dual penalty.
\begin{thm}\label{thm.3.1}
The regression-based dual penalty $\widetilde{M}(\boldsymbol{\alpha},\mathbf{z})$ defined in (\ref{eq.3.1.9}) is a feasible dual penalty, i.e., $\mathbb{E}_{0}\left[\widetilde{M}
(\boldsymbol{\alpha},\mathbf{z})\right]=0, \; \forall \boldsymbol{\alpha}\in \mathbb{A}_{\mathcal{F}}$.
\end{thm}
\begin{proof}
We have $\forall \boldsymbol{\alpha}\in \mathbb{A}_{\mathcal{F}}$
\begin{equation*}
\begin{split}
&\mathbb{E}_{0}\left[\widetilde{M}
(\boldsymbol{\alpha},\mathbf{z})\right]
=\mathbb{E}\left[\widetilde{M}(\boldsymbol{\alpha},\mathbf{z})\big|
\mathscr{F}_0\right]
=\mathbb{E}\left[\sum_{n=0}^{N-1}\sum_{i\in I} \widetilde{\beta}_{n,i}(x_n,a_n)\cdot b_i(z_{n+1})\Big|
\mathscr{F}_0\right]\\
={}&\mathbb{E}\left[\sum_{n=0}^{N-2}\sum_{i\in I} \widetilde{\beta}_{n,i}(x_n,a_n)\cdot b_i(z_{n+1})\Big|
\mathscr{F}_0\right]+\mathbb{E}\left[\sum_{i\in I} \widetilde{\beta}_{N-1,i}(x_{N-1},a_{N-1})\cdot b_i(z_{N})\Big|
\mathscr{F}_0\right]\\
\overset{(i)}={}&\mathbb{E}\left[\sum_{n=0}^{N-2}\sum_{i\in I} \widetilde{\beta}_{n,i}(x_n,a_n)\cdot b_i(z_{n+1})\Big|
\mathscr{F}_0\right]\\
&+\mathbb{E}\left[\mathbb{E}
\left[\sum_{i\in I}\widetilde{\beta}_{N-1,i}(x_{N-1},a_{N-1})\cdot b_i(z_{N})\Big|\mathscr{F}_{N-1}\right]\Big|\mathscr{F}_0\right]\\
\overset{(ii)}={}&\mathbb{E}\left[\sum_{n=0}^{N-2}\sum_{i\in I} \widetilde{\beta}_{n,i}(x_n,a_n)\cdot b_i(z_{n+1})\Big|
\mathscr{F}_0\right]+\mathbb{E}\left[\sum_{i\in I} \widetilde{\beta}_{N-1,i}(x_{N-1},
a_{N-1})\cdot\mathbb{E}\left[ b_i(z_{N})\big|\mathscr{F}_{N-1}\right]\Big|\mathscr{F}_0\right]\\
\overset{(iii)}={}&\mathbb{E}\left[\sum_{n=0}^{N-2}\sum_{i\in I} \widetilde{\beta}_{n,i}(x_n,a_n)\cdot b_i(z_{n+1})\Big|
\mathscr{F}_0\right]+\mathbb{E}\left[\sum_{i\in I} \widetilde{\beta}_{N-1,i}(x_{N-1},
a_{N-1})\cdot\mathbb{E}\left[ b_i(z_{N})\right]\Big|\mathscr{F}_0\right]\\
\overset{(iv)}={}&\mathbb{E}\left[\sum_{n=0}^{N-2}\sum_{i\in I} \widetilde{\beta}_{n,i}(x_n,a_n)\cdot b_i(z_{n+1})\Big|
\mathscr{F}_0\right]+0\\
\cdots&\\
={}&\mathbb{E}\left[
\sum_{i\in I} \widetilde{\beta}_{0,i}(x_0,a_0)\cdot b_i(z_{1})\Big|\mathscr{F}_0\right]
=\sum_{i\in I} \widetilde{\beta}_{0,i}(x_0,
a_0)\cdot\mathbb{E}\left[b_i(z_{1})\big|\mathscr{F}_0\right]\\
={}&\sum_{i\in I} \widetilde{\beta}_{0,i}(x_0,
a_0)\cdot\mathbb{E}\left[b_i(z_{1})\right]
=0,
\end{split}
\end{equation*}
where equality $(i)$ follows from the tower property of conditional expectations and $\mathscr{F}_0\subseteq \mathscr{F}_{N-1}$, equality $(ii)$ follows from the fact that $\boldsymbol{\alpha}$ is $\mathcal{F}$-adapted, and thus $\{\widetilde{\beta}_{N-1,i}(x_{N-1},
a_{N-1})\}$ are $\mathscr{F}_{N-1}$-measurable, equality $(iii)$ follows from the fact that $z_{N}$ is independent of $\mathscr{F}_{N-1}$, and finally equality $(iv)$ follows from (\ref{eq.3.3}).
\end{proof}
\noindent From the proof, we can see the fact that each basis function $b_i(\cdot)$ has expectation zero w.r.t. the probability measure $\rho$ is essential for $\widetilde{M}(\boldsymbol{\alpha},\mathbf{z})$ to maintain the dual feasibility. Theorem \ref{thm.3.1} implies that the dual problem (\ref{eq.2.5}) with penalty $\widetilde{M}(\boldsymbol{\alpha},\mathbf{z})$ provides a valid upper bound on the optimal value function, i.e., $V^{\widetilde{M}}_0(x_0)\ge V_0(x_0), \forall x_0 \in \mathscr{X}_0$. The tightness of the upper bound $V^{\widetilde{M}}_0(x_0)$ is directly affected by the accuracy of the penalty approximation, and it depends on the choice of $\mathcal{B}$ as well as the accuracy of the regression.
An ideal choice of $\mathcal{B}$ should possess the following properties: completeness, orthogonality, and countability. The reason is that completeness guarantees (\ref{eq.3.1.3}) to hold, orthogonality significantly simplifies the linear system (\ref{eq.3.1.4}) and facilitates computation of the coordinates, and countability determines whether the linear system is solvable or not. Moreover, one would prefer a basis $\mathcal{B}$ that satisfies (\ref{eq.3.3}) so that the resulted regression-based dual penalty is feasible, and thus induces valid upper bounds.
\\
\\
Following these guidelines, we consider two choices of $\mathcal{B}$ as follows. The first one is the orthonormal basis of the Hilbert $L^2$ space induced by the probability measure $\rho$. We will show that it possesses good properties of completeness, orthogonality, and countability. The second one is the basis consisting of all the centralized moments of the noise distribution. We will show that it results in a simple linear system of equations for determining the coordinates. Finally, we will show that both functional bases result in regression-based penalties that are dual feasible and good approximations of the optimal dual penalty, and thus lead to valid and tight upper bounds numerically.
\subsection{Regression Approach with $L^2$ Orthonormal Basis}\label{sec3.b:orthonormal}
To ease the presentation, let us first lay out some preliminaries for $L^2$ space. Recall the random noises $\{z_n\}$ are i.i.d. with probability measure $\rho$ on support space $\Xi \subseteq \mathbb{R}^d$. Consider the $L^2$ space induced by the measure  $\rho$, denoted by $L^2(\rho)$. In particular, $L^2(\rho)$ is the space consisting of all $\rho$-measurable scalar function $f$ such that $|f|^2$ is integrable. That is, $L^2(\rho)=\{f: \int_{\Xi} |f|^2 d\rho<\infty\}$.
\\
\\
By Theorem 4.1.3 in \cite{bogachev2007measure}, $L^2(\rho)$ is a Banach (hence complete) space equipped with the $L^2$-norm $\norm{\cdot}_2$, where
\begin{equation*}
\norm{f}_2=\sqrt{\int_{\Xi}|f|^2 d\rho},\quad \forall f\in L^2(\rho).
\end{equation*}
Here note that a metric space $S$ is called complete if every Cauchy sequence of points in $S$ has a limit in $S$. A nice property of $L^2(\rho)$ is that the $L^2$ norm $\norm{\cdot}_2$ is generated by the inner product defined by
$$
\langle f, g\rangle \overset{\triangle}= \int_{\Xi} f\cdot g d \rho,\quad \forall f, g \in L^2(\rho).
$$
Obviously, $\norm{f}_2=\sqrt{\langle f, f\rangle}$. Therefore, $L^2(\rho)$ is a Hilbert space (see, e.g., page 255 in \cite{bogachev2007measure}) equipped with the inner product defined above. Now let us introduce the definition of \emph{orthonormal basis} in a Hilbert space.
\begin{defn}\label{def.3.1}(\emph{\textbf{Orthonormal Basis}})
For a Hilbert space $H$ equipped with inner product $\langle\cdot,\cdot\rangle_H$ and norm $\norm{\cdot}_H$, a family of mutually orthogonal unit vectors $\mathcal{E}=\{e_i \in H: i\in I\}$ is called an (complete) orthonormal basis of $H$ if
\begin{itemize}
\item[(i)] $\norm{e_i}_H=1, \;\forall i \in I$; $\langle e_i, e_j\rangle_H=0,\;\forall i, j\in I, i\neq j$.
\item[(ii)] Every element $f\in H$ can be represented as the linear combination of $\{e_i\}$, i.e.,
    $$
    f\overset{H}=\sum_{i\in I} \beta_i e_i, \;\mbox{or equivalently,}\;
    \norm{f-\sum_{i\in I} \beta_i e_i}_H=0,
    $$
    where at most countably many coefficients $\beta_i$ may be nonzero.
\end{itemize}
\end{defn}
\noindent It immediately follows that if $f\overset{H}=\sum_i \beta_i e_i$, then
$$
\beta_i=\langle f, e_i\rangle_H\quad\mbox{and}\quad \norm{f}_H
=\sqrt{\sum_i |\beta_i|^2}.
$$
For a general Hilbert space $H$, by Zorn's lemma (see, e.g., \cite{bogachev2007measure}), there always exists an orthonormal basis; however, the number of basis functions in that orthonormal basis might be uncountable. Fortunately, for the specific Hilbert space $L^2(\rho)$, we have the following proposition on the countability of its orthonormal basis.
\begin{prop}\label{prop.3.1}
Let $\mathcal{E}=\{e_i\in L^2(\rho): i\in I\}$ be an orthonormal basis of $L^2(\rho)$, then $\mathcal{E}$ is countable.
\end{prop}
\begin{proof}
Note that a metric space $S$ is called separable if it contains a countable and everywhere dense subset. By Corollary 4.2.2 in \cite{bogachev2007measure}, $L^2(\rho)$ is a separable Hilbert space. Therefore, by Corollary 4.3.4 in \cite{bogachev2007measure}, any orthonormal basis of $L^2(\rho)$ is countable.
\end{proof}
\noindent Proposition \ref{prop.3.1} implies that we can rewrite the orthonormal basis $\mathcal{E}$ as $\mathcal{E}=\{e_i\in L^2(\rho): i=0,1,2,...,\}$. Without loss of generality, let us assume $e_0\equiv1$; otherwise we could simply perform the Gram–--Schmidt process to achieve that. Note that the basis $\mathcal{E}$ is uniquely determined when $e_0\equiv1$. To gain more intuition, below are examples of orthonormal basis of $L^2(\rho)$ for several common probability measures.
\\
\\
\textbf{Example 1. Finite Discrete Distribution.} Assume $\rho$ is a finite discrete probability measure with positive probability masses on $\Xi=\{y_0,...,y_{p-1}\}$. It is easy to verify that the cardinality of any orthonormal basis of $L^2(\rho)$ is $p$. In fact, $$
\mathcal{G}=\{g_i: g_i(y)=\frac{1}{\sqrt{\rho(y_i)}} \mathds{1}\{y=y_i\},\; i=0,1,...,p-1.\}
$$
is an orthonormal basis of $L^2(\rho)$.
\\
\\
\textbf{Example 2. Standard Normal Distribution.} Assume $\rho$ is the measure induced by a standard normal distribution such that $\rho(z)=1/\sqrt{2\pi}\exp(-\frac{z^2}{2}), \forall z\in \Xi=\mathbb{R}$, then the Hermite polynomials
$$
\{1;\;z;\;\frac{1}{\sqrt{2}}\left(z^2-1\right);\; \frac{1}{\sqrt{6}}\left(z^3-3z\right);...\}
$$
is a countable orthonormal basis of $L^2(\rho)$.
%
\\
\\
Next, we will show that the $L^2$ orthonormal basis $\mathcal{E}_1:=\mathcal{E}\setminus\{e_0\}$ is a valid choice of basis $\mathcal{B}$ in the proposed regression approach.
In view of the state dynamics (\ref{eq.2.1}), for a fixed $(x_n,a_n)$ pair, $V_{n+1}(x_{n+1})=V_{n+1}(f(x_n, a_n, z_{n+1}))$
is a random variable as a function of the random variable $z_{n+1}$. Further suppose $V_{n+1}(x_{n+1})$ has finite second moment, i.e.,
$V_{n+1}(f(x_n, a_n, \cdot))\in L^2(\rho)$. Since $\mathcal{E}=\{e_i: i=0,1,...\}$ with $e_0\equiv1$ is an orthonormal basis of $L^2(\rho)$, we can express $V_{n+1}(x_{n+1})$ w.r.t. $\mathcal{E}$, i.e.,
\begin{equation}\label{eq.3.10}
V_{n+1}(x_{n+1})=\sum_{i=0}^{\infty}\delta_{n,i}(x_n, a_n)\cdot e_i(z_{n+1}),
\end{equation}
where, following Definition \ref{def.3.1}, $\{\delta_{n,i}(x_n, a_n)\}$ are coordinates such that
\begin{equation}\label{eq.3.11}
\begin{split}
&\delta_{n,i}(x_n, a_n) =\langle V_{n+1}(x_{n+1}), e_i(z_{n+1})\rangle
=\int_{\Xi} V_{n+1}(x_{n+1})\cdot e_i(z_{n+1})\rho(dz_{n+1})\\
={}&\mathbb{E}\left[V_{n+1}(x_{n+1})\cdot e_i(z_{n+1})|x_n, a_n\right], \; i=0, 1,...
\end{split}
\end{equation}
In particular, noting that $e_0\equiv1$, then
\begin{equation}\label{eq.3.12}
\delta_{n,0}(x_n, a_n)=\mathbb{E}\left[V_{n+1}(x_{n+1})|x_n, a_n\right].
\end{equation}
Combining \eqref{eq.3.1.2}, \eqref{eq.3.10} and \eqref{eq.3.12}, we have
\begin{equation*}
\begin{split}
&M^{\ast}_n(\boldsymbol{\alpha},\mathbf{z})=
\left(V_{n+1}(x_{n+1})-
\mathbb{E}\left[V_{n+1}(x_{n+1})|x_n, a_n\right]\right)\\
={}&\sum_{i=0}^{\infty}\delta_{n,i}(x_n, a_n)\cdot e_i(z_{n+1})-\delta_{n,0}(x_n, a_n)\cdot e_0(z_{n+1})\\
={}&\sum_{i=1}^{\infty}\delta_{n,i}(x_n, a_n)\cdot e_i(z_{n+1}).
\end{split}
\end{equation*}
Thus, the optimal dual penalty could be rewritten as
\begin{equation}\label{eq.3.14}
M^{\ast}(\boldsymbol{\alpha},\mathbf{z})=\sum_{n=0}^{N-1}
\sum_{i=1}^{\infty}\delta_{n,i}(x_n, a_n)\cdot e_i(z_{n+1}),\quad
\mbox{where}\; a_n=\alpha_n(x_n),
\end{equation}
where $\delta_{n,i}(x_n, a_n)$ is given by (\ref{eq.3.11}). Comparing (\ref{eq.3.1.6}) with (\ref{eq.3.14}), we can see that the $L^2$ orthonormal basis $\mathcal{E}_1=\{e_1,e_2,...\}$ is a valid choice of $\mathcal{B}$ in the regression approach. Moreover, the function $h_i(\cdot)$ in (\ref{eq.3.1.5}) reduces to $e_i(\cdot)$ due to $\mathcal{E}_1$'s orthogonality. Hence, the corresponding coordinates $\{\delta_{n,i}(x_n, a_n)\}$ are straightforward to compute as in (\ref{eq.3.11}).
\\
\\
Finally, let us show that using $\mathcal{E}_1$ in the regression approach induces feasible dual penalties, and thus valid upper bounds on the optimal value function. To be more specific, let
\begin{equation*}
\widetilde{M}^l(\boldsymbol{\alpha},\mathbf{z})\overset{\triangle}=
\sum_{n=0}^{N-1}\sum_{i=1}^P \widetilde{\delta}_{n,i}(x_n,
a_n)\cdot e_i(z_{n+1}),\quad
\mbox{where}\; a_n=\alpha_n(x_n),
\end{equation*}
denote the regression-based dual penalty w.r.t. the basis $\mathcal{E}_1$, where $P$ is the order of the basis truncation, and $\{\widetilde{\delta}_{n,i}(x_n, a_n)\}$ is the estimation of $\{\delta_{n,i}(x_n, a_n)\}$ after regression. Notice that $\forall n=0,...,N-1, i=1,2,...$
\begin{equation}\label{eq.3.16}
\mathbb{E}\left[e_i(z_{n+1})\right]=\int_{\Xi}e_i(z)\rho(dz)
=\int_{\Xi}e_i(z)\cdot e_0(z)\rho(dz)=
\langle e_i(z_{n+1}), e_0(z_{n+1})\rangle=0,
\end{equation}
where we use the fact $e_0 (z_{n+1})\equiv 1$ and $\mathcal{E}$'s orthogonality. Then $\mathcal{E}_1$ is a choice of $\mathcal{B}$ that satisfies (\ref{eq.3.3}). By Theorem \ref{thm.3.1}, $\widetilde{M}^l(\boldsymbol{\alpha},\mathbf{z})$ is dual feasible.
\\
\\
As mentioned previously, the tightness of the upper bound $V^{\widetilde{M}^l}_0(x_0)$ is directly affected the accuracy of the approximation $\widetilde{M}^l(\boldsymbol{\alpha},\mathbf{z})$. If zero regression error is assumed, then $\widetilde{M}^l(\boldsymbol{\alpha},\mathbf{z})$ converges to $M^{\ast}(\boldsymbol{\alpha},\mathbf{z})$ in $L^2$ as the order of basis truncation, $P$, goes to infinity, due to
the completeness of $\mathcal{E}$. When regression error is taken into account, the construction of regressors $\boldsymbol{\phi}$ is very important in controlling that error. We will discuss this issue later.
\\
\\
Although, in theory, the unique $L^2$ orthonormal basis $\mathcal{E}_1$ could be calculated sequentially via Gram-Schmidt process for any probability measure $\rho$. In practice, this procedure might be complicated for complex probability measures, which is the drawback of using $L^2$ orthonormal basis in the regression approach. Next, we will derive an alternative functional basis that is easy to calculate. It is derived by carrying out Taylor series expansion on the value function, and consists of centralized moments of the noise distribution along each dimension.

\subsection{Regression Approach with Taylor Series Basis}\label{sec3.c:taylor}
Now let us derive the functional basis of $\mathbb{M}_\mathcal{F}$ induced by Taylor series expansion on the value function. To ease the presentation, let us further assume the state dynamics is linear in noise as follows.
\begin{equation}\label{eq.3.17}
x_{n+1}=b(x_n, a_n)+\sigma(x_n, a_n)z_{n+1},\; n=0,1,...,N-1,
\end{equation}
where $b(\cdot,\cdot)$ and $\sigma(\cdot,\cdot)$ are deterministic $\mathscr{F}$-measurable functions, and $z_{n+1}$ has finite moments (up to some order). Although the derivation could be easily extended to the setting with general state dynamics (\ref{eq.2.1}), we confine it within (\ref{eq.3.17}) to better demonstrate that several existing approaches could be regarded as special cases of the proposed approach. Note that the state dynamics (\ref{eq.3.17}) covers a wide range of dynamic programming problems. For example, it has been widely used when modeling dynamics of asset/stock prices in finance literature, where, for instance, $x$ might represent the asset/stock prices, $b(\cdot,\cdot)$ and $\sigma(\cdot,\cdot)$ might represent the drifts and the volatilities in the asset prices, respectively. The reward functions and the objective remain the same. Therefore, the primal problem is solved as before. We will focus on deriving a functional basis of the dual penalty space that facilitates approximation of the optimal dual penalty.
\\
\\
For convenience, let $\widehat{x}_{n+1}:=b(x_n, a_n)$, i.e., $\widehat{x}_{n+1}$ represents the expected future state at period $(n+1)$ given state $x_n$ and action $a_n$ at period $n$. In particular, $\widehat{x}_{n+1}$ is $\mathscr{F}_n$-measurable. Rewrite the state dynamics (\ref{eq.3.17}) as
$$
x_{n+1}-\widehat{x}_{n+1}=\sigma(x_n, a_n) z_{n+1},\; n=0,1,...,N-1.
$$
Consider the Taylor series expansion on $V_{n+1}(x_{n+1})$ around the point $\widehat{x}_{n+1}$ up to $R$-th order
\begin{equation}\label{eq.3.18}
V_{n+1}(x_{n+1})\approx V_{n+1}(\widehat{x}_{n+1})
+\sum_{r=1}^R\frac{1}{r!}\frac{\partial^r}{\partial x^r} V_{n+1}(\widehat{x}_{n+1})\cdot (\sigma)^r (z_{n+1})^r,
\end{equation}
where we use $\sigma$ as an abbreviation for $\sigma(x_n,a_n)$.
Taking expectations w.r.t. $\mathscr{F}_n$ on both sides of (\ref{eq.3.18}), we have
\begin{equation}\label{eq.3.19}
\mathbb{E}\left[V_{n+1}(x_{n+1})|x_n, a_n\right]\approx V_{n+1}(\widehat{x}_{n+1})+\sum_{r=1}^R
\frac{1}{r!}\frac{\partial^r}{\partial x^r} V_{n+1}(\widehat{x}_{n+1})\cdot (\sigma)^r \mathbb{E}\left[(z_{n+1})^r\right],
\end{equation}
where we use the fact that $V_{n+1}(\widehat{x}_{n+1})$, its partial derivatives and $\sigma(x_n,a_n)$ are $\mathscr{F}_n$-measurable, and $z_{n+1}$ is independent of $\mathscr{F}_n$. Subtracting (\ref{eq.3.19}) from (\ref{eq.3.18}), we have
\begin{equation*}
V_{n+1}(x_{n+1})-\mathbb{E}\left[V_{n+1}(x_{n+1})|x_n, a_n\right]
\approx\sum_{r=1}^R\frac{1}{r!}\frac{\partial^r}{\partial x^r} V_{n+1}(\widehat{x}_{n+1})\cdot (\sigma)^r \Bigl((z_{n+1})^r-\mathbb{E}\bigl[(z_{n+1})^r\bigr]\Bigr).
\end{equation*}
The above approximation has a nice structure because the expectations are on $(z_{n+1})^r$ instead of $V_{n+1}(x_{n+1})$, which could be calculated analytically. However, the partial derivatives could be difficult to compute directly and often require approximation methods such as finite difference. Nevertheless, it inspires us to consider the functional basis
$$
\mathcal{D}\overset{\triangle}=\{d_r: d_r(z)=(z^r-\mathbb{E}[z^r]),\; r=1,...,R\}
$$
of the dual penalty space. It follows that the optimal dual penalty $M^{\ast}(\boldsymbol{\alpha},\mathbf{z})$ could be expressed w.r.t. $\mathcal{D}$ as
\begin{equation*}
M^{\ast}(\boldsymbol{\alpha},\mathbf{z})\approx
\sum_{n=0}^{N-1}\sum_{r=1}^R\gamma_{n,r}(x_n, a_n)\cdot
\Bigl((z_{n+1})^r-\mathbb{E}\bigl[(z_{n+1})^r\bigr]\Bigr)
=\sum_{n=0}^{N-1}\sum_{r=1}^R\gamma_{n,r}(x_n, a_n)\cdot d_r(z_{n+1}),
\end{equation*}
where $\{\gamma_{n,r}(x_n, a_n)\}$ are the corresponding coordinates. Following the derivations in (\ref{eq.3.1.3})-(\ref{eq.3.1.4}), we obtain a simple linear system of equations for $\{\gamma_{n,r}(x_n, a_n)\}$ as follows:
\begin{equation}\label{eq.3.20}
\begin{bmatrix}
\mathbb{E}_n\bigl[V_{n+1}(x_{n+1})\bigl (z_{n+1}-\mathbb{E}[z_{n+1}]\bigr)\bigr]\\
\mathbb{E}_n\bigl[V_{n+1}(x_{n+1})\bigl ((z_{n+1})^2-\mathbb{E}\left[(z_{n+1})^2\right]\bigr)\bigr]\\
\vdots\\
\mathbb{E}_n\bigl[V_{n+1}(x_{n+1})\bigl ((z_{n+1})^R-\mathbb{E}\left[(z_{n+1})^R\right]\bigr)\bigr]
\end{bmatrix}
\approx
\begin{bmatrix}
Cov\begin{pmatrix}
z_{n+1}\\
(z_{n+1})^2\\
\vdots\\
(z_{n+1})^R
\end{pmatrix}
\end{bmatrix}
\begin{bmatrix}
\gamma_{n,1}(x_n, a_n)\\
\gamma_{n,2}(x_n, a_n)\\
\vdots\\
\gamma_{n,R}(x_n, a_n)
\end{bmatrix},
\end{equation}
where $Cov(\cdot)$ represents the covariance matrix of a random vector. Note that it is nonsingular, and thus we can solve (\ref{eq.3.20}) as
\begin{equation}\label{eq.3.21}
\begin{bmatrix}
\gamma_{n,1}(x_n, a_n)\\
\gamma_{n,2}(x_n, a_n)\\
\vdots\\
\gamma_{n,R}(x_n, a_n)
\end{bmatrix}
\approx
\begin{bmatrix}
Cov\begin{pmatrix}
z_{n+1}\\
(z_{n+1})^2\\
\vdots\\
(z_{n+1})^R
\end{pmatrix}
\end{bmatrix}^{\mathbf{-1}}
\begin{bmatrix}
\mathbb{E}_n\bigl[V_{n+1}(x_{n+1})\bigl (z_{n+1}-\mathbb{E}[z_{n+1}]\bigr)\bigr]\\
\mathbb{E}_n\bigl[V_{n+1}(x_{n+1})\bigl ((z_{n+1})^2-\mathbb{E}\left[(z_{n+1})^2\right]\bigr)\bigr]\\
\vdots\\
\mathbb{E}_n\bigl[V_{n+1}(x_{n+1})\bigl ((z_{n+1})^R-\mathbb{E}\left[(z_{n+1})^R\right]\bigr)\bigr]
\end{bmatrix}.
\end{equation}
For simplicity, let us rewrite (\ref{eq.3.21}) explicitly as
\begin{equation}\label{eq.3.22}
\gamma_{n,r}(x_n, a_n)\approx\mathbb{E}
\left[V_{n+1}(x_{n+1})l_r(z_{n+1})|x_n, a_n\right],
\quad r=1,...,R,
\end{equation}
where $\{l_r(\cdot)\}$ are simple polynomial functions of degrees $\le R$ because the vector on the right hand side of (\ref{eq.3.21}) only involves monomials of $z_{n+1}$ of degrees $\le R$. Therefore, the coordinates $\{\gamma_{n,r}(x_n, a_n)\}$ could be easily computed, which is the main advantage of using the function basis $\mathcal{D}$ in the regression approach. The drawback is that it does not guarantee convergence of the optimal dual penalty approximation as $R$ goes to infinity, because in general Taylor series expansion is not convergent.
\\
\\
Finally, similar to using the function basis $\mathcal{E}_1$, let us show that using $\mathcal{D}$ in the regression approach produces feasible dual penalties, and thus valid upper bounds on the optimal value function. To this end, let
\begin{equation*}
\widetilde{M}^t(\boldsymbol{\alpha},\mathbf{z})\overset{\triangle}=
\sum_{n=0}^{N-1}\sum_{r=1}^R \widetilde{\gamma}_{n,i}(x_n,
a_n)\cdot d_i(z_{n+1}),
\quad \mbox{where}\; a_n=\alpha_n(x_n)
\end{equation*}
denote the regression-based dual penalty generated via Algorithm \ref{alg.3.1} with $\mathcal{D}$ as the functional basis $\mathcal{B}$, where $\{\widetilde{\gamma}_{n,i}(x_n, a_n)\}$ is the estimation of $\{\gamma_{n,i}(x_n, a_n)\}$ through regression. Note that for $n=0,...,N-1, r=1,...,R$,
\begin{equation*}
\mathbb{E}[d_r(z_{n+1})]=
\mathbb{E}[(z_{n+1})^R-\mathbb{E}[(z_{n+1})^R]]
=\mathbb{E}[(z_{n+1})^R]-\mathbb{E}[(z_{n+1})^R]=0.
\end{equation*}
Thus, $\mathcal{D}$ is a choice of $\mathcal{B}$ that satisfies (\ref{eq.3.3}). Therefore, by Theorem \ref{thm.3.1} $\widetilde{M}^t(\boldsymbol{\alpha},\mathbf{z})$ is a feasible dual penalty. Similarly, an immediate implication is that the dual problem (\ref{eq.2.5}) with penalty $\widetilde{M}^t(\boldsymbol{\alpha},\mathbf{z})$ provides a valid upper bound on the optimal value function, i.e., $V^{\widetilde{M}^t}_0(x_0)\ge V_0(x_0), \forall x_0 \in \mathscr{X}_0$. Again, the tightness of the upper bound $V^{\widetilde{M}^t}_0(x_0)$ is directly affected by the accuracy of the approximation $\widetilde{M}^t(\boldsymbol{\alpha},\mathbf{z})$, which heavily relies on the regressors $\boldsymbol{\phi}$ constructed in the regression. We will discuss this issue next.
\begin{rem}\label{rem.3.1}
Interestingly, one could verify that if the noises $\{z_n\}$ are symmetric, then the first two basis functions $d_1, d_2\in \mathcal{D}$ coincide with the first two basis functions $e_1, e_2\in\mathcal{E}_1$ (after taking out the normalization factors), respectively. For example, when $\{z_n\}$ follows a standard normal distribution, then $\{d_1=z, d_2=(z^2-1)\}\subset \mathcal{D}$ coincides with $\{e_1=z, e_2=\frac{1}{\sqrt{2}}(z^2-1)\}\subset \mathcal{E}_1$ in Example 2 of Section \ref{sec3.b:orthonormal}.
\end{rem}

\subsection{Implementation}
We will discuss two important implementation issues regarding the proposed regression approach, i.e., how to properly construct the regressors in the regression and how to solve the inner optimization problem in the resulted dual problem.
\\
\\
As mentioned previously, the accuracy of the regression-based dual penalty $\widetilde{M}(\boldsymbol{\alpha},\mathbf{z})$ heavily relies on the choice of regressors $\boldsymbol{\phi}$ for regressing the coordinates $\beta$. The general guideline is to construct $\boldsymbol{\phi}$ similar to $\beta$ in structure. In view of the equation $\beta_{n,i}(x_n, a_n)=\mathbb{E}\left[V_{n+1}(x_{n+1})\cdot h_i(z_{n+1})|x_n, a_n\right]$, a natural candidate regressor is $\mathbb{E}\left[\widetilde{V}_{n+1}(x_{n+1})\cdot h_i(z_{n+1})|x_n, a_n\right]$, where $\widetilde{V}_{n+1}(\cdot)$ is an approximation of $V_{n+1}(\cdot)$ with a closed-form expression. For example, $\widetilde{V}_{n+1}(\cdot)$ could be the value function induced by a naive policy, or an approximate optimal value function, or the value function of a simpler DP (e.g., through a relaxation of the constraints on the actions of the original DP). Multiple regressors of such form could be included in the regression. More specifically, if the Taylor series basis $\mathcal{D}$ is used, then there exists good regressors of a simpler form. Note that the coordinates $\{\gamma_{n,r}(x_n, a_n)\}$ take the positions of the partial derivatives in the Taylor series expansion of the optimal dual penalty. Hence, a good candidate regressor is $\frac{\partial^r}{\partial x^r} \widetilde{V}_{n+1}(\widehat{x}_{n+1})\cdot (\sigma)^r$, where $\widetilde{V}_{n+1}(\cdot)$ admits partial derivatives that are easy to derive. Again, $\widetilde{V}_{n+1}(\cdot)$ could be the value function induced by a naive policy, or an approximate optimal value function, or the value function of a simpler DP.
\\
\\
Theoretically, after plugging the regression-based dual penalty $\widetilde{M}(\boldsymbol{\alpha},\mathbf{z})$ into the dual problem (\ref{eq.2.5}), we only need to solve a series of deterministic inner optimization problem to generate an upper bound estimator. However, they might be hard to solve without good structural properties such as convexity. Therefore, for a convex DP, ideally the dual penalty plugged in the dual problem should preserve convexity in the inner optimization problem. Note that the optimal dual penalty $M^{\ast}(\boldsymbol{\alpha},\mathbf{z})$ might not be convex even if the value functions $\{V_n(x_n)\}$ are concave, and in this case plugging in the optimal dual penalty causes the inner optimization problem to lose convexity. One solution to this issue is to linearize the optimal dual penalty around a fixed policy. The resulted optimal dual penalty approximation is affine in policy, and thus preserves convexity in the dual problem.
\\
\\
The proposed framework of regression approach does not suffer from this issue. It is robust in the sense that an optimal dual penalty approximation with desired structural properties could be generated by using proper regressors, noting that the resulted dual penalty is a linear combination of the regressors. For example, for a convex DP, the regression-based dual penalty preserves convexity in the dual problem by constructing regressors that are affine in policy.

\section{Special Cases}\label{sec4:special}
It turns out, in a broader sense, several existing approaches to approximating the optimal dual penalty are special cases of the proposed regression approach under specific settings.

\subsection{American Option Pricing}\label{sec4.a:option}
Optimal stopping problem is one of the simplest non-trivial types of DPs. Various methodologies and numerical methods are first proposed under the setting of optimal stopping, and then generalized to general DPs. American option pricing, as one of the most salient applications of optimal stopping, has been extensively studied. We will demonstrate that the non-nested simulation approach developed in \cite{Belomestny:2009} to approximating the optimal dual martingale in American option pricing could be regarded as a special case of the proposed regression approach with the first-order Taylor series basis $\{d_1\}$.
\\
\\
The problem setting is as follows. Assume the asset price $X(t)$ satisfies a stochastic differential equation (SDE) w.r.t. Brownian motion
\begin{equation*}
dX\left( t \right) = b\left( {t,X\left( t \right)} \right)dt + \sigma \left( {t,X\left( t \right)} \right)dW\left( t \right),
\end{equation*}
where $t\in[0,T]$, $X(t)$ denotes the asset price at time $t$ with given initial deterministic value $X(0)=X_{0}$, $W(t)$ represents the (standard) Brownian motion, and the coefficients $b$, $\sigma$ are functionals satisfying mild regularity conditions. By convention, we use $\{\mathscr{F}_{t}:0\leq t \leq T\}$ to denote the augmented information filtration generated by the Brownian motion $\{W(t)\}$. Consider a Bermudan option on $X(t)$ that can be exercised at any date from the time set $\mathcal{T}  = \{ {T_0},{T_1},...,{T_N }\}$, with $0=T_{0}<T_{1}<\cdot\cdot\cdot<T_{N}=T$. When exercising at time $T_{n}\in \mathcal{T}$, the option holder receives a payoff
${H_{{T_n}}}= h( {{T_n},X( {{T_n}} )} )$. The goal is to evaluate the price of the Bermudan option, that is, to find
\begin{equation*}
(P):\;V_{0}=\sup_{\tau \in \mathcal{T}}\mathbb{E}\left[h\left(\tau,X\left(\tau \right)\right)|X\left(0\right)=X_{0}\right],
\end{equation*}
where $\tau$ is an exercising strategy, i.e., a stopping time adapted to the filtration $\mathcal{F}=\{\mathscr{F}_{T_{n}}:n=0,...,N\}$ taking values in $\mathcal{T}$, and $V_{0}$ denotes the Bermudan option price at time $T_{0}$ with initial asset price $X_{0}$.
\\
\\
As for the dual formulation, \cite{andersen:2004} and \cite{haugh:2004} show that all the $\mathcal{F}$-adapted martingales in $\mathbb{M}_{\mathcal{F}}=\{M=(M_0,...,M_N): M\;\text{is $\mathcal{F}$-adapted}\}$ are feasible dual penalties. In particular, the optimal dual martingale $M^\ast$ is the Doob-Meyer martingale component of the Bermudan price process $\{V_{T_n}: n=0,...,N\}$, i.e.,
\begin{equation}\label{eq.4.3}
\left\{\begin{array}{l}
M^\ast_{0}=0,
\\
M^\ast_{T_{n+1}}=M^\ast_{T_{n}}+V_{T_{n+1}}-\mathbb{E}_{T_{n}}
\left[V_{T_{n+1}}\right], n=0,...,N-1.
\end{array}\right.
\end{equation}
The question is how to efficiently estimate $M^\ast$ in a non-nested manner. Note that the martingale $M^\ast$ is driven by Brownian motion, by martingale representation theorem we have
\begin{equation*}
M^\ast_{T_{n}} = \int_{0}^{T_{n}} {{C_s}dW\left(s\right)}, \; n=0,...,N,
\end{equation*}
where ${C_s}$ is a predictable process w.r.t. the filtration $\mathcal{F}$. \cite{Belomestny:2009} proposes an Ito sum approximation scheme as follows.
\begin{equation*}
M^\ast_{T_{n}} = \int_{0}^{T_{n}} {C_s}dW(s)\approx \sum\limits_{i=0}^{n-1} C_{T_i}\left(\Delta W_{T_{i+1}}\right), n=0,...,N,
\end{equation*}
where $\left(\Delta W_{T_{i+1}}\right):=\left(W(T_{i + 1})-W(T_{i})\right)$ denotes the Brownian increment. Combining with (\ref{eq.4.3}), we have
\begin{equation}\label{eq.4.6}
V_{T_{n+1}}-\mathbb{E}_{T_{n}}\left[V_{T_{n+1}}\right]= M^\ast_{T_{n+1}}-M^\ast_{T_{{n}}} 
\approx  C_{T_n}\left(\Delta W_{T_{n+1}}\right),~\; n=0,...,N-1.
\end{equation}
Multiplying both sides of \eqref{eq.4.6} by the Brownian increment $\left(\Delta W_{T_{n+1}}\right)$ and taking conditional expectations w.r.t. $\mathscr{F}_{T_{n}}$, we obtain
\begin{equation*}
C_{T_{n}} \approx \frac{1}{T_{n+1}-T_{n}}\mathbb{E}_{T_{n}}
\left[V_{T_{n+1}}\left(\Delta W_{T_{n+1}}\right)\right],
~~n=0,....,N-1.
\end{equation*}
Starting from here, regression is further applied to estimate $C_{T_{n}}$ in order to avoid nested simulation. Eventually, a regression-based dual martingale $\widetilde{M}=(\widetilde{M}_{T_0}, \widetilde{M}_{T_1},...,
\widetilde{M}_{T_N})$ is given by
\begin{equation}\label{eq.4.8}
\widetilde{M}_{T_n}\overset{\triangle}=
\sum_{i=0}^{n-1}\widetilde{C}_{T_i}\cdot
\left(\Delta W_{T_{i+1}}\right),\quad n=0,...,N,
\end{equation}
where $\widetilde{C}_{T_i}$ is the estimation of $C_{T_i}$ after regression. Since the celebrated Ito's Lemma and Martingale Representation Theorem could be viewed as the results of carrying out Taylor series expansion on a stochastic function, an Ito sum approximation of a stochastic integral could be viewed as a first-order Taylor series expansion type of scheme. Thus, approximation scheme (\ref{eq.4.8}) could be viewed as a special case of the proposed regression approach with the first-order Taylor series basis $\{d_1\}$.
\begin{rem}
\cite{Belomestny:2009} actually considers a finer Ito sum approximation of the optimal dual martingale by introducing a finer partition of the time span. Here we present a simpler version of their approach so that its connection with our proposed approach is more clear.
\end{rem}

\subsection{Controlled Markov Diffusion}\label{sec4.b:diffusion}
\cite{ye2015information} studies the form of optimal dual penalty under the setting of controlled Markov diffusion (CMD) and show that it is a stochastic integral. It then inspires the authors to propose an approximation scheme for the optimal dual penalty of a discrete-time DP with state dynamics (\ref{eq.3.17}) as follows:
\begin{equation}\label{eq.4.9}
\widehat{M}(\boldsymbol{\alpha},\mathbf{z})\overset{\triangle}=
\sum_{n=0}^{N-1}\frac{\partial}{\partial x} V_n(x_n)\cdot\sigma(x_n, a_n)\cdot z_{n+1},
\end{equation}
where the partial derivative $\frac{\partial}{\partial x} V_n(x_n)$ is estimated via finite difference method. This scheme is derived by mimicking the Ito sum approximation of the optimal dual penalty (which is a stochastic integral) under the setting of CMD. The connection with our proposed regression approach is evident in the sense that the approximation scheme (\ref{eq.4.9}) could be interpreted as a result of a first-order Taylor series expansion of the optimal dual penalty, and the coordinates are estimated using finite difference method instead of regression. Another difference we must point out is that $\frac{\partial}{\partial x} V_{n+1}(\hat{x}_{n+1})$ instead of $\frac{\partial}{\partial x} V_n(x_n)$ is used in our proposed regression approach. We will show that, at least for the following linear quadratic control (LQC) problem, our approach is better in terms of approximation accuracy.
\subsection{Linear Quadratic Control}\label{sec4.c:linear}
Now let us consider the Taylor series expansion of the optimal dual penalty for the classic LQC problem, since the LQC problem has been extensively studied and applied, see \cite{bertsekas:2007}. Assuming a linear state dynamics
\begin{equation*}
x_{n+1}=A_{n}x_{n}+B_{n}a_{n}+z_{n+1},\;~ n=0,...,N-1,
\end{equation*}
with an expected quadratic cost
\begin{equation}\label{eq.4.11}
\mathbb{E}_{0}\left[x_{N}^{T}Q_{N}x_{N}
+\sum\limits_{n=0}^{N-1}\left(x_{n}^{T}Q_{n}x_{n}+
a_{n}^{T}R_{n}a_{n}\right)\right],
\end{equation}
where $A_{n}, B_{n}, Q_{n}, R_{n}$ are known matrices of appropriate dimensions, $\{z_{n}\}$ are proper-dimensional i.i.d. vector random variables with zero mean and finite second moment, and $Q_{n}$ and $R_{n}$ are positive semi-definite. The objective is to select a non-anticipative policy $\boldsymbol{\alpha}^{\ast}$ that minimizes the total costs defined in (\ref{eq.4.11}). It turns out that the optimal control policy (e.g., see \cite{bertsekas:2007}) $\alpha^{\ast}_{n}$ admits a closed-form expression $\alpha^{\ast}_{n}(x_n)=L_{n}x_{n}$, where the gain matrix $L_{n}$ is given by
$L_{n}=-\left(B_{n}^{T}K_{n+1}B_{n}+
R_{n}\right)^{-1}B_{n}^{T}K_{n+1}A_{n}$,
and the positive semi-definite matrix $K_{n}$ is given recursively by
\begin{equation*}
\left\{\begin{array}{l}
K_{N}=Q_{N},
\\
K_{n}=A^{T}_{n}\Bigl(K_{n+1}-K_{n+1}B_{n}
\left(B_{n}^{T}K_{n+1}B_{n}+R_{n}\right)^{-1}
B_{n}^{T}K_{n+1}\Bigr)A_{n}+Q_{n}.
\end{array}\right.
\end{equation*}
The optimal value function $V_{n}(x_n)$ then takes the form
\begin{equation*}
V_{n}(x_{n})=x^{T}_{n}K_{n}x_{n}+
\sum\limits_{i=n}^{N-1}\mathbb{E}
\left[z^{T}_{n+1}K_{n+1}z_{n+1}\right],\; n=0,...,N.
\end{equation*}
For the corresponding dual problem via information relaxation, it is easy to verify that (see, e.g., \cite{haugh2012linear}) the value-based optimal dual penalty $M^{\ast}(\boldsymbol{\alpha},\mathbf{z})$ takes the form
\begin{flalign*}
M^{\ast}(\boldsymbol{\alpha},\mathbf{z})&=\sum\limits_{n=0}^{N-1}
\Bigl(V_{n+1}(x_{n+1})-\mathbb{E}
\left[V_{n+1}(x_{n+1})|x_n, a_n\right]\Bigr)\nonumber\\
&=\sum\limits_{n=0}^{N-1}\Bigl\{
2(A_{n}x_{n}+B_{n}a_{n})^{T}K_{n+1}z_{n+1}
+z^{T}_{n+1}K_{n+1}z_{n+1}-\mathbb{E}
\left[z^{T}_{n+1}K_{n+1}z_{n+1}\right]\Bigr\}.
\end{flalign*}
Now let us approximate $M^{\ast}(\boldsymbol{\alpha},\mathbf{z})$ using the second-order Taylor series basis $\{d_1, d_2\}$, where we apply second-order Taylor series expansion on $V_{n+1}(x_{n+1})$ around $\widehat{x}_{n+1}=A_{n}x_{n}+B_{n}a_{n}$.  That is,
\begin{eqnarray*}
&&\widetilde{M}(\boldsymbol{\alpha},\mathbf{z})\\
&=&\sum\limits_{n=0}^{N-1}
\Biggl\{\frac{\partial}{\partial x}V^T_{n+1}(\widehat{x}_{n+1})  z_{n+1}+z_{n+1}^T \frac{\partial^2}{\partial x^2}V_{n+1}(\widehat{x}_{n+1})z_{n+1}-\mathbb{E}\left[z_{n+1}^T \frac{\partial^2}{\partial x^2}V_{n+1}(\widehat{x}_{n+1})z_{n+1}\right]\Biggr\}\nonumber\\
&=&\sum\limits_{n=0}^{N-1}\Bigl\{
2(A_{n}x_{n}+B_{n}a_{n})^{T}K_{n+1}z_{n+1}
+z^{T}_{n+1}K_{n+1}z_{n+1}-\mathbb{E}
\left[z^{T}_{n+1}K_{n+1}z_{n+1}\right]\Bigr\}
=M^{\ast}(\boldsymbol{\alpha},\mathbf{z}).
\end{eqnarray*}
Therefore, the second-order Taylor series expansion of the optimal dual penalty is exact, which is due to the linearity of state dynamics and quadratic structure of the optimal value function. If using the approximation scheme in (\ref{eq.4.9}), we would have
\begin{equation*}
\widehat{M}(\boldsymbol{\alpha},\mathbf{z})
=\sum_{n=0}^{N-1}\left(\frac{\partial}{\partial x} V_n(x_n) \right)^T z_{n+1}
=\sum_{n=0}^{N-1}2x_n^T K_n z_{n+1}\neq M^{\ast}(\boldsymbol{\alpha},\mathbf{z}).
\end{equation*}

\section{Numerical Experiments: Dynamic Trading}\label{sec5:numerical}

\subsection{Dynamic Trading with Predictable Returns and Transaction Costs}\label{sec5.a:trading}
To demonstrate the effectiveness of the proposed regression approach in generating accurate approximations of optimal dual penalty and tight upper bounds, we will study the dynamic trading problem with predictable returns and transaction costs in \cite{garleanu2013dynamic}. This model has been empirically tested in real financial markets and exhibits nice structural properties. The problem formulation is as follows.
\\
\\
Consider an investor that attempts to trade $D$ securities at each time $t\in\{0,1,...,T\}$. The securities' price changes (returns) are driven by $K$ market factors, i.e.,
\begin{equation*}
r_{t+1}=\mu_t+B f_t + z^{(1)}_{t+1}, ~~t=0,...,T-1,
\end{equation*}
where $r_{t+1}$ is the $D\times 1$ vector of security returns at time $t+1$, $\mu_t$ is the deterministic ``risk-free'' return, $f_t$ is the $K\times 1$ vector of market factors that predict returns, $B$ is the $D\times K$ matrix of constant factor loadings, and $\{z^{(1)}_{t+1}\}$ are i.i.d. zero-mean $D\times 1$ random vectors with covariance matrix $Var(z^{(1)}_{t+1})=\Sigma$ that represent the unpredictable noise in the return. Further assume the market return vector $f_t$ is already known to the investor at time $t$ (before the trading) and it follows a self-evolving state dynamics
\begin{equation}\label{eq.5.2}
\Delta f_t \overset{\triangle}=f_{t+1}-f_{t}=-\Phi f_t + z^{(2)}_{t+1},~~ t=0,...,T-1,
\end{equation}
where $\Phi$ is a $K\times K$ matrix of mean-reversion coefficients for the factors, and $\{z^{(2)}_{t+1}\}$ are i.i.d. $K\times 1$ zero-mean random vectors with covariance matrix $Var(z^{(2)}_{t+1})=\Psi$ that represent the shock affecting the predictors. Usually, we assume $\{z^{(2)}_{t+1}\}$ are normally distributed. Assume $\Psi$ satisfies the standard conditions such that $f$ is stationary. Note that the state dynamics (\ref{eq.5.2}) could be rewritten as $f_{t+1}=(I-\Phi)f_t+z^{(2)}_{t+1}$, where $I$ is the $K\times K$ identity matrix.
\\
\\
Assuming trading is costly in the sense that it impacts the market and results in a price move. In particular, \cite{garleanu2013dynamic} argues that the transaction cost associated with trading quantity $a_{t}=x_t-x_{t-1}$, where $x_t$ is the share quantity held by the investor at time $t$ (after trading) and $a_t$ is the trading quantity at time $t$, is given by $TC(a_t)=\frac{1}{2}a_t^T \Lambda a_t$. Here $\Lambda$ is a $D\times D$ deterministic positive semi-definite matrix measuring the level of trading cost. The intuition is that trading $a_t$ shares will move the price by $\frac{1}{2}\Lambda a_t$ and result in a total trading cost $(\frac{1}{2}\Lambda a_t)^T a_t$, which is exactly $TC(a_t)$.
\\
\\
The investor's objective is to choose a (non-anticipative) trading strategy that maximizes the total expected excess return with the trading costs and risk penalties taken into account, i.e.,
\begin{equation}\label{eq.5.3}
(P):\;V_1(x_0, f_1)=\max_{\boldsymbol{\alpha}\in \mathbb{A}_\mathcal{F}}
\mathbb{E}_1\biggl[\sum_{t=1}^T\Bigl(x_t^T B f_t -\frac{1}{2}a_t^T \Lambda a_t-\frac{\gamma}{2}x_t^T \Sigma x_t\Bigr)\biggr],
\end{equation}
where $\mathbb{A}_\mathcal{F}$ describes the set of all admissible policies that satisfy certain trading constraints such as no short-sells or complete liquidation of the initial position by the end of trading horizon $T$, $\mathcal{F}=\{\mathscr{F}_0,...,\mathscr{F}_T\}$ is the natural information filtration, $x_t^T B f_t$ is the expected intermediate excess return, and $\gamma$ is the risk-aversion coefficient.
\\
\\
Notice that the primal problem (\ref{eq.5.3}) falls into the realm of LQC if no trading constraints are imposed (except for complete liquidation of the initial position). Therefore, both the optimal value function and the optimal trading strategy could be solved in closed form using backward dynamic programming. For example, the optimal value function is quadratic in the state $(x_{t-1}, f_t)$ and the optimal control is an affine function of the state. In particular, for risk-neutral ($\gamma=0$) trading, \cite{moallemi2015dynamic} shows that the optimal value function $J_{t+1}(x_{t},f_{t+1})$ satisfies
\begin{equation}\label{eq.5.4}
\left\{\begin{array}{l}
J_{T}(x_{T-1},f_T)=-\frac{1}{2}x_{T-1}^T\Lambda x_{T-1}, \\
J_{t}(x_{t-1},f_{t})=-\frac{1}{2}x_{t-1}^T A_{xx,t} x_{t-1} +
x_{t-1}^T A_{xf,t}f_{t}
+\frac{1}{2}f_{t}^T A_{ff,t}f_{t}+A_{t},\;t=T-1,...,1,
\end{array}\right.
\end{equation}
where $A_{xx,t}$, $A_{xf,t}$, $A_{ff,t}$ and $A_{t}$ are coefficient matrices of proper dimensions that satisfy the following backward recursions:
\begin{align*}
&A_{xx,T}=\Lambda, \quad A_{xf,T}=\mathbf{0},\quad A_{ff,T}=\mathbf{0},\quad A_{T}=0;\\
&A_{xx,t}=-\Lambda(\Lambda+A_{xx,t+1})^{-1}\Lambda+\Lambda,\\
&A_{xf,t}=\Lambda(\Lambda+A_{xx,t+1})^{-1}(B+A_{xf,t+1}(I-\Phi)),\\
&A_{ff,t}=(B+A_{xf,t+1}(I-\Phi))^T(\Lambda+A_{xx,t+1})^{-1}
(B+A_{xf,t+1}(I-\Phi))+(I-\Phi)^T A_{ff,t+1}(I-\Phi),\\
&A_{t}=trace(\Psi A_{ff,t+1})+A_{t+1}.
\end{align*}
Furthermore, the optimal policy $\boldsymbol{\alpha}^\ast$ is given by
\begin{equation}\label{eq.5.5}
\alpha^{\ast}_t(x_{t-1})=(\Lambda+A_{xx,t+1})^{-1}\bigl(\Lambda x_{t-1}+(B+A_{xf,t+1}(I-\Phi))f_t\bigr)-x_{t-1}.
\end{equation}
Note that after imposing trading constraints, as many real-world trading problems do, the classic theory in LQC fails and in general it is difficult to find an optimal policy.
Here let us consider a very natural set of trading constraints: only selling of securities is allowed, short-sell of the securities is not allowed, and the initial position must be completely liquidated by terminal horizon. Mathematically, the constraint set $\mathscr{A}$ could be formulated as
\begin{equation*}
\mathscr{A}\overset{\triangle}=\bigl\{\mathbf{a}=(a_1,...,a_T): x_t=x_{t-1}+a_t, a_t\le 0, x_t\ge 0, t=1,...,T; x_T=0\bigr\}
\end{equation*}
Notice that all the constraints in $\mathscr{A}$ are linear (hence convex) in actions. Thus, the problem maintains to be tractable when deterministic convex objective functions are of interest.
\\
\\
In the existence of trading constraints $\mathscr{A}$, the optimal trading strategy is not readily available. Alternatively, we focus on finding sufficiently good suboptimal strategies. Although the optimal policy of the unconstrained problem is unlikely to be feasible for the constrained problem (since one of the constraints in $\mathscr{A}$ might be violated), a natural feasible policy could be derived by projecting the optimal policy of the unconstrained problem onto the constraint set $\mathscr{A}$. In particular, we consider the following projected linear quadratic control (PLQC) policy $\widehat{\boldsymbol{\alpha}}=(\widehat{\alpha}_1,...,
\widehat{\alpha}_T)$ with
\begin{equation}\label{eq.5.7}
\widehat{\alpha}_t (x_{t-1})=\max\{-x_{t-1},\; \min\{0,\;\alpha^{\ast}_t(x_{t-1})\}\},~~ t=1,...,T,
\end{equation}
where the projection is componentwise and $\alpha^{\ast}_t(x_{t-1})$ is given by (\ref{eq.5.5}). From (\ref{eq.5.7}), the projection of control onto $\mathscr{A}$ is easy to implement, and thus it is straightforward to evaluate the PLQC policy via Monte Carlo simulation. Eventually, it results in a lower bound $V^{\widehat{\boldsymbol{\alpha}}}_1(x_0, f_1)$ on the optimal value function $V_1(x_0, f_1)$.
\begin{rem}\label{rem.4.1}
In \cite{moallemi2015dynamic}, the authors also consider a time-weighted average price (TWAP) policy where an equal quantity of shares ($a_t=-x_0/T$) is traded at every period, a deterministic (DETER) policy where a deterministic trading strategy is given by solving the deterministic trading problem with all noises ignored, a model predictive control (MPC) policy where the trading strategy is given by solving a sequence of deterministic trading problems without noises, and a linear rebalancing (LRB) policy where the control is given by the best rebalancing policy that is affine in all the available market predictors. Our preliminary numerical tests show that these three policies (TWAP, DETER, MPC) perform significantly worse than the PLQC policy. Although the LRB policy is competitive (often worse) compared with the PLQC policy in terms of performance; it is more sensitive to formulation, more complex in structure, and more difficult to solve. Therefore, we focus on the easy-to-implement PLQC policy.
\end{rem}

\subsection{Dual Formulation with Regression-based Penalties}\label{sec5.b:dual}
Now let us consider the dual formulation of primal problem (\ref{eq.5.3}). Following the derivation in Section \ref{sec2.a: dual formulation}, we know that the dual problem
\begin{equation}\label{eq.5.8}
(D):\; V^M_1(x_0, f_1)=
\mathbb{E}_1\Biggl[\max_{\boldsymbol{\alpha}\in \mathbb{A}}
\sum_{t=1}^T\biggl(x_t^T B f_t -\frac{1}{2}a_t^T \Lambda a_t
-\frac{\gamma}{2}x_t^T \Sigma x_t\biggr)-M(\boldsymbol{\alpha},\mathbf{z}^{(2)})\Biggr],
\end{equation}
where $\mathbf{z}^{(2)}:=(z_{1}^{(2)},...,z_{T}^{(2)})$ and $M(\boldsymbol{\alpha},\mathbf{z}^{(2)})$ is a feasible dual penalty. It follows that $V^M_1(x_0, f_1)\ge V_1(x_0, f_1)$. In particular, when plugging in the optimal dual penalty $M^{\ast}(\boldsymbol{\alpha},\mathbf{z}^{(2)})$ the exact value function is recovered, i.e., $V^{M^{\ast}}_1(x_0, f_1)=V_1(x_0, f_1)$, where
\begin{equation*}
M^{\ast}(\boldsymbol{\alpha},\mathbf{z}^{(2)})
=\sum\limits_{t=1}^{T-1}\Bigl(V_{t+1}(x_{t}, f_{t+1})-
\mathbb{E}\left[V_{t+1}(x_{t}, f_{t+1})|x_{t-1}, a_{t}, f_{t}\right]\Bigr).
\end{equation*}
Since the optimal value functions $\{V_{t+1}(x_{t}, f_{t+1})\}$ are not available, we have to approximate $M^{\ast}(\boldsymbol{\alpha},\mathbf{z}^{(2)})$. In particular, we will apply the proposed regression approach. Note that the objective function is quadratic, we will use the second-order Taylor series basis (which is equivalent to the second-order $L^2$ orthonormal basis for this problem) in the algorithm.
\\
\\
Following the derivation in Section \ref{sec3.c:taylor}, let us express the optimal dual penalty w.r.t. the basis $\{d_1=z, d_2=z^2-\mathbb{E}[z^2]\}$, i.e.,
\begin{equation*}
M^{\ast}(\boldsymbol{\alpha},\mathbf{z}^{(2)})\approx
\sum_{t=1}^{T-1}\biggl(
\gamma^T_{t,1}(x_{t-1}, a_{t}, f_{t}) z^{(2)}_{t+1}
+\gamma^T_{t,2}(x_{t-1}, a_{t}, f_{t})
\left((z^{(2)}_{t+1})^2-
\mathbb{E}\left[(z^{(2)}_{t+1})^2\right]\right)\biggr),
\end{equation*}
where $\{\gamma_{t,r}(x_{t-1}, a_{t}, f_{t})\}$ are computed via
\begin{equation}\label{eq.5.11}
\left\{\begin{array}{l}
\gamma_{t,1}(x_{t-1}, a_{t}, f_{t})=Cov(z^{(2)}_{t+1})^{-1}\mathbb{E}_t\left[
V_{t+1}(x_{t}, f_{t+1})z^{(2)}_{t+1} \right]\\
\gamma_{t,2}(x_{t-1}, a_{t}, f_{t})=Cov((z^{(2)}_{t+1})^2)^{-1}\mathbb{E}_t\Bigl[
V_{t+1}(x_{t}, f_{t+1})
\left((z^{(2)}_{t+1})^2
-\mathbb{E}\left[(z^{(2)}_{t+1})^2\right]\right) \Bigr].
\end{array}\right.
\end{equation}
Note that we are slightly abusing the notation $(z^{(2)}_{t+1})^2$ to denote the componentwise square of a random vector. We will use Algorithm \ref{alg.3.1} to estimate $\{\gamma_{t,r}(x_{t-1}, a_{t}, f_{t})\}$ in (\ref{eq.5.11}). The key is how to construct proper regressors in the regression. As discussed in the previous section, a good candidate regressor is $\frac{\partial^r}{\partial f^r} \widetilde{V}_{t+1}(x_{t},\widehat{f}_{t+1})\cdot (\sigma)^r$, where $\widetilde{V}_{t+1}(\cdot)$ is either the value function of a naive policy or the value function of a simpler DP such that its derivatives admit closed-form expressions. Here $\widehat{f}_{t+1}=(I-\Phi)f_{t}$. Recall that the corresponding unconstrained problem is an LQC problem, and its value functions $\{J_{t+1}(x_t, f_{t+1})\}$ admit closed-form expressions as in (\ref{eq.5.4}). Therefore, we can use $\frac{\partial^r}{\partial f^r} J_{t+1}(x_{t},\widehat{f}_{t+1})\cdot (\sigma)^r$ as one of the regressors for regressing $\gamma_{n,r}(x_{t-1}, a_{t}, f_{t})$. Moreover, $J_{t+1}(x_{t},\widehat{f}_{t+1})$ is quadratic in $(x_{t},\widehat{f}_{t+1})$, $x_{t}$ is affine in $(x_{t-1}, a_{t})$, and $\widehat{f}_{t+1}$ is linear in $f_{t}$. Thus, the first-order and second-order derivatives of $J_{t+1}(x_{t},\widehat{f}_{t+1})$ w.r.t. $f$ are both affine in $(x_{t-1}, a_{t}, f_{t})$. It follows that the regressor, and thus the resulted dual penalty is affine in policy and preserves convexity in the dual problem. In particular, our numerical tests show that the regressors
\begin{equation*}
\boldsymbol{\phi}_{t,1}(x_{t-1}, a_{t}, f_{t})=\Bigl\{1, \frac{\partial}{\partial f} J_{t+1}(x_{t},\widehat{f}_{t+1})=
A^T_{xf,t+1}(x_{t-1} + a_{t})+A_{ff,t+1}(I-\Phi)f_{t}\Bigr\}
\end{equation*}
and
\begin{equation}\label{eq.5.13}
\boldsymbol{\phi}_{t,2}=\left\{\frac{\partial^2}{\partial f^2} J_{t+1}(x_{t},\widehat{f}_{t+1})=
diag(A_{ff,t+1})\right\}
\end{equation}
are sufficiently good for regressing $\gamma_{n,1}(x_{t-1}, a_{t}, f_{t})$ and $\gamma_{n,2}(x_{t-1}, a_{t}, f_{t})$, respectively. Here note that the derivatives are componentwise, $diag(A_{ff,t+1})$ denotes the diagonal vector of matrix $A_{ff,t+1}$, and we suppress the dependence of function $\boldsymbol{\phi}_{t,2}$ on $(x_{t-1}, a_{t}, f_{t})$ because $diag(A_{ff,t+1})$ is a constant.  To this end, let us denote the resulted regression-based dual penalty as
\begin{equation}\label{eq.5.14}
\widetilde{M}^t(\boldsymbol{\alpha},\mathbf{z}^{(2)})
=\sum_{t=1}^{T-1}\Bigl(
\widetilde{\gamma}^T_{t,1}(x_{t-1}, a_{t}, f_{t})  z^{(2)}_{t+1}
\widetilde{\gamma}^T_{t,2}(x_{t-1}, a_{t}, f_{t})
\left((z^{(2)}_{t+1})^2-
\mathbb{E}\left[(z^{(2)}_{t+1})^2\right]\right)\Bigr),
\end{equation}
where $\widetilde{\gamma}_{t,r}(x_{t-1}, a_{t}, f_{t})$ is the estimation of $\gamma_{t,r}(x_{t-1}, a_{t}, f_{t})$ after regression. Notice that the regressor $\boldsymbol{\phi}_{t,2}$ in (\ref{eq.5.13}) does not depend on $(x_{t-1}, a_{t}, f_{t})$, it follows that $\widetilde{\gamma}^T_{t,2}$ does not depend on $(x_{t-1}, a_{t},f_{t})$. Therefore, the second component of $\widetilde{M}^t(\boldsymbol{\alpha},\mathbf{z}^{(2)})$,
$$
\sum_{t=1}^{T-1}
\widetilde{\gamma}^T_{t,2}
\left((z^{(2)}_{t+1})^2-
\mathbb{E}\left[(z^{(2)}_{t+1})^2\right]\right)
$$
does not involve states and actions. Thus, it does not affect the optimization in the inner optimization problems. In fact, it plays the role of control variate in the sense that it has zero mean and thus does not affect the expectation of the upper bound estimator; however, the correlation between this term and the reward function might help reduce its variance. Plugging $\widetilde{M}^t(\boldsymbol{\alpha},\mathbf{z}^{(2)})$ into the dual problem (\ref{eq.5.8}), we have
\begin{equation*}
V^{\widetilde{M}^t}_1(x_0, f_1)=
\mathbb{E}_1\biggl[\max_{\boldsymbol{\alpha}\in \mathbb{A}}
\sum_{t=1}^T\Bigl(x_t^T B f_t -\frac{1}{2}a_t^T \Lambda a_t
-\frac{\gamma}{2}x_t^T \Sigma x_t\Bigr)-\widetilde{M}^t(\boldsymbol{\alpha},\mathbf{z}^{(2)})\biggr].
\end{equation*}
Note that $\widetilde{M}^t(\boldsymbol{\alpha},\mathbf{z}^{(2)})$ is a feasible dual penalty, and thus $V^{\widetilde{M}^t}_1(x_0, f_1)\ge V_1(x_0, f_1)$. To obtain an unbiased estimator of $V_1^{\widetilde{M}^t}(x_0, f_1)$, we can simulate $L$ i.i.d. sample paths of $\mathbf{z}^{(2)}$, solve the deterministic inner optimization problem corresponding to every sample path, and finally take the sample average of the optimal values. Since the inner optimization problems are deterministic and convex, we use the CVX package (Mosek solver) in Matlab to solve them.

\subsection{Numerical Results}\label{sec5.c:results}
We will use the model parameters calibrated in \cite{moallemi2015dynamic}. However, the authors only consider the case of trading one stock (Apple, Inc), i.e., $D=1$. To exhibit the effectiveness of our regression approach for high-dimensional DPs, we ``replicate'' their model parameters to obtain a high-dimensional dynamic trading problem. The model parameters are summarized as follows.
\\
\\
Time horizon is $T=12$ or $24$. The number of stocks to be liquidated is $D=5, 10$, or $25$. The initial position of every stock to be liquidated is $[x_0^1,...,x_0^D]$, a vector in $\mathbb{R}^D$ with $x_0^d=10000$. There are $K=2$ market factors to predict the return on each stock, each with a different mean reversion speed. The return of the stocks follows the dynamics
\begin{equation*}
r_{t+1}^d=0.0726+0.3375f_{t,1}-0.0720f_{t,2}+z^{(1)}_{t+1,d},~~
d=1,...,D;\; t=0,...,T-1,
\end{equation*}
where $\Sigma=Var(z^{(1)}_{t+1,d})=0.048$. Thus, the $D\times K$ matrix $$B=[0.3375, -0.072; 0.3375, -0.072;...,0.3375, -0.072].$$
Similarly, the market factor $f$ follows the dynamics
\begin{equation*}
\left\{\begin{array}{l}
f_{t+1,1}-f_{t,1}=-0.5f_{t,1}+z^{(2)}_{t+1,1},\\
f_{t+1,2}-f_{t,2}=-0.7f_{t,2}+z^{(2)}_{t+1,2},
\end{array}\right.
\end{equation*}
where
$$\Psi=Var(z^{(2)}_{t+1})= \begin{bmatrix}0.0379 &0\\
0&0.0947
\end{bmatrix},$$
and
$$ \Phi=\begin{bmatrix} 0.5 & 0\\
0&0.7
\end{bmatrix}.
$$
We also test other choices of $\Phi$ in the appendix.
Instead of generating $f_0, f_1$ randomly, we let $f_0=[1,1]'$ and $f_1=(I-\Phi)f_0$. The transaction cost matrix $\Lambda=\lambda\cdot\widetilde{\Lambda}$, where $\lambda=2.14\times 10^{-5}$, and $\widetilde{\Lambda}=\Gamma\Gamma^T$. Here $\Gamma$ is an upper triangular matrix such that
$$
\Gamma=\begin{bmatrix}
\frac{1}{\sqrt{D}}&\frac{1}{\sqrt{D}}&\cdots &\frac{1}{\sqrt{D}}\\
0 & \frac{1}{\sqrt{D-1}} & \cdots &\frac{1}{\sqrt{D-1}}\\
\vdots & \vdots &\ddots &\vdots\\
0 & 0& \cdots & 1
\end{bmatrix}_{D\times D}.
$$
The reason for using such a matrix is to ensure that diagonal elements of $\widetilde{\Lambda}$ are all one, meaning the transaction costs are the same across all the stocks. We also test other choices of $\lambda$ in the appendix.
\\
\\
With the above model parameters, we solve for the parameters (i.e., $A_{xx}$, $A_{xf}$, $A_{ff}$, and $A$) of the optimal policy to the corresponding unconstrained problem. Then the PLQC policy is computed via policy projection (\ref{eq.5.5}). To evaluate the PLQC policy and generate a lower bound estimator, we run a simulation of $M=10^6$ i.i.d. sample paths with the same initial state $(x_0, f_1)$, exercise the PLQC policy along each sample path, and compute the corresponding accumulated reward. The value function of the PLQC policy, i.e., the lower bound on the optimal value function, is estimated by taking the average of those reward samples. As for the upper bound estimation, we consider the upper bounds induced by three different feasible dual penalties: (1) Zero penalty (note that it is a trivial feasible dual penalty), meaning the investor has a perfect foresight of all the future information without any penalization; (2) Regression-based penalty with the first-order Taylor series basis, i.e.,
\begin{equation}\label{eq.5.16}
\widetilde{M}^t_1(\boldsymbol{\alpha},\mathbf{z}^{(2)})=
\sum_{t=1}^{T-1}
\widetilde{\gamma}^T_{t,1}(x_{t-1}, a_{t}, f_{t}) z^{(2)}_{t+1};
\end{equation}
(3) Regression-based penalty with the second-order Taylor series basis, i.e., $\widetilde{M}^t(\boldsymbol{\alpha},\mathbf{z}^{(2)})$ in  (\ref{eq.5.14}).
To regress $\{\gamma^T_{t,r}(x_{t-1}, a_{t}, f_{t})\}$, we implement Algorithm \ref{alg.3.1} with the PLQC policy and the $M=10^6$ sample paths generated in estimating the lower bound. After constructing all three dual penalties, we simulate $L=400$ sample paths to generate the upper bound estimators. $L$ is much smaller than $M$ because the variances of the upper bound estimators induced by good approximations of the optimal dual penalty are small.
To illustrate the performance of the PLQC policy, we compute its duality gap, which is the ratio of the difference between the lower bound and the tightest upper bound to the lower bound. If the duality gap is small enough, we could claim the policy is sufficiently good. The detailed results are summarized in Table \ref{table.5.1}.
\begin{table}[!htb]
\centering
\caption{ Dual Bounds of PLQC Policy with Parameter $\Phi$. } \label{table.5.1}
\begin{tabular}{|c|c|c|c|c|c|c|} \hline
$D$& $T$ & LB & UB 1 & UB 2 & UB 3& DualGap\\
& & PLQC & ZP & F-O  & S-O &\\
& & $(\$k)$ & $(\$k)$ & $(\$k)$ & $(\$k)$ & \\
\hline\hline
$5$&$12$& $14.937$ & $18.096$ & $15.365$ & $15.263$ & $2.18\%$ \\
       &    & $(0.074)$ & $(1.304)$ & $(0.107)$& $(0.055)$&    \\
\hline
$10$&$12$& $24.303$ & $30.240$ & $24.544$ & $24.645$ & $0.99\%$\\
       &    & $(0.140)$ & $(2.224)$ & $(0.152)$& $(0.070)$&    \\
\hline
$25$&$12$& $32.090$ & $44.124$ & $32.657$ & $32.529$ & $1.37\%$\\
       &    & $(0.324)$ & $(5.237)$ & $(0.211)$& $(0.109)$&    \\
\hline
$5$&$24$& $18.635$ & $28.060$ & $20.012$ & $20.103$ & $7.39\%$\\
       &    & $(0.103)$ & $(2.094)$ & $(0.250)$& $(0.152)$&    \\
\hline
$10$&$24$& $33.309$ & $47.738$ & $35.099$ & $35.285$ & $5.37\%$\\
       &    & $(0.212)$ & $(3.816)$ & $(0.306)$& $(0.193)$&    \\
\hline
$25$&$24$& $62.971$ & $86.537$ & $64.401$ & $65.224$ & $2.27\%$\\
       &    & $(0.528)$ & $(9.673)$ & $(0.370)$& $(0.201)$&    \\
\hline
\end{tabular}
\\
{\small The experiment parameters: $M=10^6$ sample paths are simulated to evaluate the lower bounds and compute the coefficients in the regression; $L=400$ sample paths are simulated to evaluate the upper bounds.}
\end{table}

\noindent In Table \ref{table.5.1}, Column 1 records the number of securities $D$, including the cases where $D=5, 10$ or $25$. Column 2 records the number of trading horizons $T$, including the cases $T=12$ or $24$. Column 3 records the lower bounds (LB) induced by the PLQC policy. Column 4 records the upper bounds (UB 1) induced by the zero penalty. Column 5 records the upper bounds (UB 2) induced by the dual penalty $\widetilde{M}^t_1(\boldsymbol{\alpha},\mathbf{z}^{(2)})$ in (\ref{eq.5.16}). Column 6 records the upper bounds (UB 3) induced by the dual penalty $\widetilde{M}^t(\boldsymbol{\alpha},\mathbf{z}^{(2)})$ in (\ref{eq.5.14}). The half confidence interval widths of the lower bound and each upper bound are presented in the parentheses. The last column records the best duality gaps (DualGap) achieved by comparing the lower bounds and the tightest upper bounds. We have the following observations:
\begin{itemize}
\item In general, in view of the small duality gaps in most cases, we could claim that the simple PLQC policy is a sufficiently good policy.
\item Comparing UB 1, UB 2 and UB 3 induced by the zero penalty, the first-order regression-based penalty, and the second-order regression-based penalty, respectively, we notice: 1) UB 1 has large duality gaps, meaning the zero penalty performs poorly as one would expect. 2) For most parameter setups, UB 2 and UB 3 have very small duality gaps, meaning that the policy used in the regression (here is the PLQC policy) is near optimal and the regression-based dual penalties are accurate approximations of the optimal dual penalty. 3) It is hard to distinguish between UB 2 and UB 3 in terms of tightness. The reason is that the inner optimization problems are the same, except for a control variate term in the second-order regression-based penalty. This term helps reduce the variance of the upper bound estimator, which is verified by the fact that the half confidence interval widths of UB 3 are narrower than the ones of UB 2.
\item We also observe that, as the number of trading horizons $T$ increases, the duality gaps increase. It indicates the possibility that the PLQC policy becomes less optimal for larger $T$. Another possibility is that the approximation of the optimal dual penalty is less accurate for larger $T$.
\item Lastly, as the number of securities in position increases, the duality gaps decrease instead of increasing as one would expect. One possible explanation is that the PLQC policy is closer to optimal for higher-dimensional problems due to the strong correlations across the assets in the transaction cost matrix $\Lambda$. Therefore, the optimal trading strategy of the unconstrained problem is more conservative, and thus less optimality is lost in the policy projection.
\end{itemize}
Overall, we can see that the PLQC policy is near optimal for the model parameters tested. More importantly, the regression-based dual penalties perform well. It demonstrates the effectiveness and efficiency of the proposed framework of regression approach in solving the dual problems of high-dimensional DPs.

\section{Conclusions}\label{sec6:conclusion}
In this paper, we develop a framework of regression approach to approximating the optimal dual penalty in general DPs, by studying the structure of the dual penalty space. The proposed framework circumvents the issue of nested simulation suffered by some of the existing approaches, and thus improves computational efficiency. Furthermore, the proposed framework requires minimal extra simulation and computational costs by reusing the samples generated during the estimation of lower bound. Lastly, the proposed framework is robust. Dual penalties with desired structural properties could be generated by constructing proper regressors in the regression. The application to a high-dimensional dynamic trading problem demonstrates its effectiveness in generating good feasible dual penalties and tight upper bounds on the optimal value function.

\bibliographystyle{ormsv080}
\bibliography{Zhou-Bibtex}
\newpage

\section*{Appendix}
The performances of the regression-based penalties under different choices of factor persistence $\Phi$. In particular,
$$
\Phi_1=\begin{bmatrix}
0.3 & 0\\
0 & 0.5
\end{bmatrix},
\Phi_2=\begin{bmatrix}
0.3 & 0\\
0 & 0.7
\end{bmatrix},
\Phi_3=\begin{bmatrix}
0.5 & 0\\
0 & 0.3
\end{bmatrix},
\Phi_4=\begin{bmatrix}
0.7 & 0\\
0 & 0.5
\end{bmatrix}.
$$

\begin{table}[htpb]
\centering
\begin{centering}
\caption{ Lower and Upper Bounds of PLQC policy. } \label{table.a.1}
\footnotesize{\begin{tabular}{|c|c|c|c|c|c|c|c|} \hline
$\Phi$ &$D$& $T$ & Lower Bound & UpperBound 1 & UpperBound 2 & UpperBound 3& Duality Gap\\
& & & PLQC Policy & Zero Penalty & First-Order  & Second-Order &\\
& & & $(\$k)$ & $(\$k)$ & $(\$k)$ & $(\$k)$ & \\
\hline\hline
$\Phi_1$ &$1$&$12$& $7.564$ & $8.327$ & $7.760$ &  $7.732$ & $2.22\%$ \\
& &                &$(0.026)$ &$(0.448)$&$(0.070)$&$(0.046)$&    \\

$\Phi_1$ &$5$&$12$& $33.571$ & $36.675$ & $33.692$ &  $33.670$ & $0.29\%$ \\
& &                &$(0.122)$ &$(1.890)$&$(0.191)$&$(0.103)$&    \\

$\Phi_1$ &$10$&$12$& $59.391$ & $63.451$ & $59.893$ & $59.837$ & $0.75\%$ \\
& &                &$(0.230)$ &$(3.509)$&$(0.291)$&$(0.136)$&    \\

$\Phi_1$ &$25$&$12$& $110.277$ & $131.213$ & $112.164$ & $113.531$ & $1.71\%$ \\
& &                &$(0.506)$ &$(8.102)$&$(0.552)$&$(0.326)$&    \\

\hline

$\Phi_2$ &$1$&$12$& $7.964$ & $9.311$ & $8.155$ &  $8.170$ & $2.40\%$ \\
& &                &$(0.026)$ &$(0.815)$&$(0.066)$&$(0.042)$&    \\

$\Phi_2$ &$5$&$12$& $35.460$ & $39.136$ & $35.955$ &  $35.758$ & $0.84\%$ \\
& &                &$(0.122)$ &$(1.787)$&$(0.241)$&$(0.096)$&    \\

$\Phi_2$ &$10$&$12$& $63.287$ & $70.992$ & $63.770$ & $63.885$ & $0.76\%$ \\
& &                &$(0.228)$ &$(3.658)$&$(0.261)$&$(0.117)$&    \\

$\Phi_2$ &$25$&$12$& $119.974$ &$135.117$ & $123.230$ & $122.642$ & $2.22\%$ \\
& &                &$(0.505)$ &$(9.307)$&$(0.526)$&$(0.277)$&    \\

\hline
$\Phi_3$ &$1$&$12$& $2.565$ &$3.917$ & $2.791$ & $2.849$ & $8.81\%$ \\
& &                &$(0.014)$ &$(0.222)$&$(0.063)$&$(0.049)$&    \\

$\Phi_3$ &$5$&$12$& $9.776$ &$13.097$ & $10.072$ & $10.374$ & $3.02\%$ \\
& &                &$(0.070)$ &$(1.217)$&$(0.150)$&$(0.084)$&    \\

$\Phi_3$ &$10$&$12$& $13.946$ &$20.531$ & $14.388$ & $14.489$ & $3.17\%$ \\
& &                &$(0.135)$ &$(2.373)$&$(0.193)$&$(0.107)$&    \\

$\Phi_3$ &$25$&$12$& $6.897$ &$14.942$ & $7.392$ & $7.372$ & $6.89\%$ \\
& &                &$(0.327)$ &$(5.140)$&$(0.250)$&$(0.142)$&    \\

\hline
$\Phi_4$ &$1$&$12$& $1.274$ &$2.029$ & $1.463$ & $1.445$ & $13.42\%$ \\
& &                &$(0.010)$ &$(0.166)$&$(0.040)$&$(0.026)$&    \\

$\Phi_4$ &$5$&$12$& $3.609$ &$5.978$ & $3.715$ & $3.756$ & $2.94\%$ \\
& &                &$(0.051)$ &$(0.817)$&$(0.081)$&$(0.043)$&    \\

$\Phi_4$ &$10$&$12$& $2.638$ &$8.173$ & $2.741$ & $2.707$ & $2.62\%$ \\
& &                &$(0.099)$ &$(1.657)$&$(0.162)$&$(0.064)$&    \\

$\Phi_4$ &$25$&$12$& $-18.907$ &$-5.989$ & $-18.356$ & $-18.448$ & $2.43\%$ \\
& &                &$(0.245)$ &$(3.906)$&$(0.159)$&$(0.081)$&    \\

\hline
\end{tabular}}
\end{centering}
\\
{\quad \small The experiment parameters: $M=10^6$ sample paths are simulated to evaluate the lower bounds and generate the regression-based penalties; $L=400$ sample paths are simulated to generate the upper bounds.}
\end{table}

\newpage
\noindent The performances of the regression-based penalties with different choices of transaction cost parameter $\lambda$. In particular, we let
$$
\lambda_1=1.07\times 10^{-5},\;
\lambda_2=2.67\times 10^{-5},\;
\lambda_3=3.21\times 10^{-5},\;
\lambda_4=4.28\times 10^{-5}.
$$
\begin{table}[htpb]
\centering
\caption{ Lower and Upper Bounds of PLQC policy. } \label{table.a.2}
\footnotesize{\begin{tabular}{|c|c|c|c|c|c|c|c|} \hline
$\lambda$ &$D$& $T$ & Lower Bound & UpperBound 1 & UpperBound 2 & UpperBound 3& Duality Gap\\
& & & PLQC Policy & Zero Penalty & First-Order  & Second-Order &\\
& & & $(\$k)$ & $(\$k)$ & $(\$k)$ & $(\$k)$ & \\
\hline\hline
$\lambda_1$ &$1$&$12$& $3.757$ & $4.691$ & $4.104$ &  $4.027$ & $7.19\%$ \\
& &                &$(0.015)$ &$(0.250)$&$(0.078)$&$(0.056)$&    \\

$\lambda_1$ &$5$&$12$& $16.963$ & $21.195$ & $17.210$ & $17.353$ & $1.46\%$ \\
& &                &$(0.075)$ &$(1.324)$&$(0.157)$&$(0.094)$&    \\

$\lambda_1$ &$10$&$12$& $29.899$ & $36.332$ & $30.377$ & $29.749$ & $1.60\%$ \\
& &                &$(0.145)$ &$(2.444)$&$(0.226)$&$(0.104)$&    \\

$\lambda_1$ &$25$&$12$& $55.321$ & $65.456$ & $56.405$ & $55.902$ & $1.05\%$ \\
& &                &$(0.338)$ &$(6.056)$&$(0.343)$&$(0.201)$&    \\

\hline

$\lambda_2$ &$1$&$12$& $3.583$ & $4.286$ & $3.759$ & $3.756$ & $4.83\%$ \\
& &                &$(0.015)$ &$(0.234)$&$(0.045)$&$(0.034)$&    \\

$\lambda_2$ &$5$&$12$& $14.119$ & $18.721$ & $14.437$ & $14.433$ & $2.25\%$ \\
& &                &$(0.072)$ &$(1.221)$&$(0.090)$&$(0.040)$&    \\

$\lambda_2$ &$10$&$12$& $21.854$ & $25.148$ & $22.444$ & $22.458$ & $2.70\%$ \\
& &                &$(0.136)$ &$(2.096)$&$(0.120)$&$(0.062)$&    \\

$\lambda_2$ &$25$&$12$& $21.162$ & $28.615$ & $21.743$ & $21.862$ & $2.75\%$ \\
& &                &$(0.320)$ &$(5.446)$&$(0.200)$&$(0.117)$&    \\

\hline
$\lambda_3$ &$1$&$12$& $3.533$ & $4.552$ & $3.626$ & $3.642$ & $2.63\%$ \\
& &                &$(0.015)$ &$(0.251)$&$(0.039)$&$(0.023)$&    \\

$\lambda_3$ &$5$&$12$& $13.326$ & $18.707$ & $13.694$ & $13.649$ & $2.42\%$ \\
& &                &$(0.072)$ &$(1.083)$&$(0.096)$&$(0.039)$&    \\

$\lambda_3$ &$10$&$12$& $19.759$ & $26.425$ & $20.091$ & $20.254$ & $1.68\%$ \\
& &                &$(0.134)$ &$(2.290)$&$(0.099)$&$(0.045)$&    \\

$\lambda_3$ &$25$&$12$& $10.798$ & $20.543$ & $11.275$ & $11.164$ & $3.39\%$ \\
& &                &$(0.318)$ &$(4.439)$&$(0.161)$&$(0.085)$&    \\

\hline
$\lambda_4$ &$1$&$12$& $3.417$ & $4.486$ & $3.523$ & $3.510$ & $2.72\%$ \\
& &                &$(0.015)$ &$(0.273)$&$(0.035)$&$(0.018)$&    \\

$\lambda_4$ &$5$&$12$& $12.008$ & $16.255$ & $12.275$ & $12.241$ & $1.94\%$ \\
& &                &$(0.070)$ &$(1.120)$&$(0.079)$&$(0.031)$&    \\

$\lambda_4$ &$10$&$12$& $15.811$ & $19.592$ & $16.215$ & $16.182$ & $2.35\%$ \\
& &                &$(0.132)$ &$(2.070)$&$(0.092)$&$(0.043)$&    \\

$\lambda_4$ &$25$&$12$& $-9.375$ & $3.890$ & $-9.258$ & $-9.335$ & $0.43\%$ \\
& &                &$(0.315)$ &$(4.766)$&$(0.155)$&$(0.108)$&    \\

\hline
\end{tabular}}
\\
{\quad \small The experiment parameters: $M=10^6$ sample paths are simulated to evaluate the lower bounds and generate the regression-based penalties; $L=400$ sample paths are simulated to generate the upper bounds.}
\end{table}

\end{document}